\documentclass[a4paper,12pt,leqno]{article}
\usepackage{amsmath,amssymb,latexsym,amscd,amsthm,amsxtra,graphicx}
\begin{document}




\newtheorem{thm}{Theorem}[subsection]
\newtheorem{lem}[thm]{Lemma}
\newtheorem{cor}[thm]{Corollary}
\newtheorem{prop}[thm]{Proposition}
\newtheorem{prob}[thm]{Problem}
\newtheorem{clm}{Claim}
\newtheorem{dfn}[thm]{Definition}
\newtheorem{conj}[thm]{Conjecture}
\newtheorem{wh}[thm]{Working hypothesis}


\newcommand{\pf}{{\it Proof} \,\,\,\,\,\,}
\newcommand{\remark}{{\it Remark} \,\,\,\,}

\newcommand{\eqb}{\begin{equation}}
\newcommand{\eqe}{\end{equation}}

\newcommand{\cpx}{{\mathbb C}}
\newcommand{\rel}{{\mathbb R}}
\newcommand{\rat}{{\mathbb Q}}
\newcommand{\itg}{{\mathbb Z}}
\newcommand{\nat}{{\mathbb N}}
\newcommand{\qrt}{{\mathbb H}}
\newcommand{\qtr}{{\mathbb H}}
\newcommand{\dbP}{{\mathbb P}}
\newcommand{\dbS}{{\mathbb S}}
\newcommand{\dbO}{{\mathbb O}}
\newcommand{\dbK}{{\mathbb K}}

\newcommand{\mathbbP}{{\mathbb P}}
\newcommand{\mathbbA}{{\mathbb A}}
\newcommand{\mathbbB}{{\mathbb B}}
\newcommand{\mathbbC}{{\mathbb C}}
\newcommand{\mathbbD}{{\mathbb D}}
\newcommand{\mathbbE}{{\mathbb E}}
\newcommand{\mathbbF}{{\mathbb F}}
\newcommand{\mathbbG}{{\mathbb G}}
\newcommand{\mathbbH}{{\mathbb H}}
\newcommand{\mathbbI}{{\mathbb I}}
\newcommand{\mathbbJ}{{\mathbb J}}
\newcommand{\mathbbK}{{\mathbb K}}
\newcommand{\mathbbL}{{\mathbb L}}
\newcommand{\mathbbM}{{\mathbb M}}
\newcommand{\mathbbN}{{\mathbb N}}
\newcommand{\mathbbO}{{\mathbb O}}
\newcommand{\mathbbQ}{{\mathbb Q}}
\newcommand{\mathbbR}{{\mathbb R}}
\newcommand{\mathbbS}{{\mathbb S}}
\newcommand{\mathbbT}{{\mathbb T}}
\newcommand{\mathbbU}{{\mathbb U}}
\newcommand{\mathbbV}{{\mathbb V}}
\newcommand{\mathbbW}{{\mathbb W}}
\newcommand{\mathbbX}{{\mathbb X}}
\newcommand{\mathbbY}{{\mathbb Y}}
\newcommand{\mathbbZ}{{\mathbb Z}}

\newcommand{\arn}{\Vec{n}}

\newcommand{\uDl}{\underline{\Delta}}
\newcommand{\uPi}{\underline{\Pi}}
\newcommand{\oDl}{\overline{\Delta}}
\newcommand{\uW}{\underline{W}}
\newcommand{\oW}{\overline{W}}
\newcommand{\uqca}{\underline{\cal Q}}
\newcommand{\upca}{\underline{\cal P}}
\newcommand{\oqca}{\overline{\cal Q}}
\newcommand{\opca}{\overline{\cal P}}
\newcommand{\urho}{\underline{\rho}}
\newcommand{\orho}{\overline{\rho}}
\newcommand{\ubbb}{\underline{\mathfrak b}}
\newcommand{\uqqq}{\underline{\mathfrak q}}
\newcommand{\uppp}{\underline{\mathfrak p}}
\newcommand{\uvvv}{\underline{\mathfrak v}}
\newcommand{\unnn}{\underline{\mathfrak n}}
\newcommand{\obbb}{\overline{\mathfrak b}}
\newcommand{\onnn}{\overline{\mathfrak n}}
\newcommand{\ola}{\overline{\lambda}}
\newcommand{\dola}{\overline{\overline{\lambda}}}
\newcommand{\ula}{\underline{\lambda}}
\newcommand{\dula}{\underline{\underline{\lambda}}}
\newcommand{\hla}{{\hat{\lambda}}}

\newcommand{\bGL}{{\mathbb GL}}

\newcommand{\bfa}{\mbox{\underline{\bf a}}}
\newcommand{\bfb}{\mbox{\underline{\bf b}}}
\newcommand{\bfc}{\mbox{\underline{\bf c}}}
\newcommand{\bfd}{\mbox{\underline{\bf d}}}
\newcommand{\bfe}{\mbox{\underline{\bf e}}}
\newcommand{\bff}{\mbox{\underline{\bf f}}}
\newcommand{\bfg}{\mbox{\underline{\bf g}}}
\newcommand{\bfh}{\mbox{\underline{\bf h}}}
\newcommand{\bfi}{\mbox{\underline{\bf i}}}
\newcommand{\bfj}{\mbox{\underline{\bf j}}}
\newcommand{\bfk}{\mbox{\underline{\bf k}}}
\newcommand{\bfl}{\mbox{\underline{\bf l}}}
\newcommand{\bfm}{\mbox{\underline{\bf m}}}
\newcommand{\bfn}{\mbox{\underline{\bf n}}}
\newcommand{\bfo}{\mbox{\underline{\bf o}}}
\newcommand{\bfp}{\mbox{\underline{\bf p}}}
\newcommand{\bfq}{\mbox{\underline{\bf q}}}
\newcommand{\bfr}{\mbox{\underline{\bf r}}}
\newcommand{\bfs}{\mbox{\underline{\bf s}}}
\newcommand{\bft}{\mbox{\underline{\bf t}}}
\newcommand{\bfu}{\mbox{\underline{\bf u}}}
\newcommand{\bfv}{\mbox{\underline{\bf v}}}
\newcommand{\bfw}{\mbox{\underline{\bf w}}}
\newcommand{\bfx}{\mbox{\underline{\bf x}}}
\newcommand{\bfy}{\mbox{\underline{\bf y}}}
\newcommand{\bfz}{\mbox{\underline{\bf z}}}
\newcommand{\bfmu}{{\underline{\bf \mu}}}

\newcommand{\bpa}{\mbox{\underline{\bf a}}}
\newcommand{\bpb}{\mbox{\underline{\bf b}}}
\newcommand{\bpc}{\mbox{\underline{\bf c}}}
\newcommand{\bpd}{\mbox{\underline{\bf d}}}
\newcommand{\bpe}{\mbox{\underline{\bf e}}}
\newcommand{\bpf}{\mbox{\underline{\bf f}}}
\newcommand{\bpg}{\mbox{\underline{\bf g}}}
\newcommand{\bph}{\mbox{\underline{\bf h}}}
\newcommand{\bpi}{\mbox{\underline{\bf i}}}
\newcommand{\bpj}{\mbox{\underline{\bf j}}}
\newcommand{\bpk}{\mbox{\underline{\bf k}}}
\newcommand{\bpl}{\mbox{\underline{\bf l}}}
\newcommand{\bpm}{\mbox{\underline{\bf m}}}
\newcommand{\bpn}{\mbox{\underline{\bf n}}}
\newcommand{\bpo}{\mbox{\underline{\bf o}}}
\newcommand{\bbpp}{\mbox{\underline{\bf p}}}
\newcommand{\bpq}{\mbox{\underline{\bf q}}}
\newcommand{\bpr}{{\mbox{\underline{\bf r}}}}
\newcommand{\bps}{\mbox{\underline{\bf s}}}
\newcommand{\bpt}{\mbox{\underline{\bf t}}}
\newcommand{\bpu}{\mbox{\underline{\bf u}}}
\newcommand{\bpv}{\mbox{\underline{\bf v}}}
\newcommand{\bpw}{\mbox{\underline{\bf w}}}
\newcommand{\bpx}{\mbox{\underline{\bf x}}}
\newcommand{\bpy}{\mbox{\underline{\bf y}}}
\newcommand{\bpz}{\mbox{\underline{\bf z}}}

\newcommand{\eps}{\varepsilon}

\newcommand{\T}{\Theta}
\newcommand{\Th}{\Theta}
\newcommand{\Ta}{\Theta^\alpha}

\newcommand{\zzz}{{\mathfrak z}}
\newcommand{\fff}{{\mathfrak f}}
\newcommand{\www}{{\mathfrak w}}
\newcommand{\vvv}{{\mathfrak v}}
\renewcommand{\ggg}{{\mathfrak g}}
\newcommand{\gggg}{{\mathfrak g}}
\newcommand{\kkk}{{\mathfrak k}}
\newcommand{\aaa}{{\mathfrak a}}
\newcommand{\saa}{\mbox{}^s{\mathfrak a}}
\newcommand{\aas}{{\mathfrak a}_{\T}}
\newcommand{\aasa}{{\mathfrak a}_{\T^\alpha}}
\newcommand{\bbb}{{\mathfrak b}}
\newcommand{\sbb}{\mbox{}^s{\mathfrak p}}
\newcommand{\spp}{\mbox{}^s{\mathfrak p}}
\newcommand{\llll}{{\mathfrak l}}
\newcommand{\llls}{{\mathfrak l}_{\T}}

\newcommand{\tll}{\tilde{\mathfrak l}}
\newcommand{\eee}{{\mathfrak e}}
\newcommand{\ccc}{{\mathfrak c}}
\newcommand{\ttt}{{\mathfrak t}}
\newcommand{\stt}{\mbox{}^s{\mathfrak t}}
\newcommand{\lls}{{\mathfrak l}_{\T}}
\newcommand{\obb}{\bar{\mathfrak b}}
\newcommand{\aad}{{\mathfrak a}^\ast_{\T}}
\newcommand{\ads}{{\mathfrak a}^\ast_{\T}}
\newcommand{\nnn}{{\mathfrak n}}
\newcommand{\nns}{{\mathfrak n}_{\T}}
\newcommand{\uuu}{{\mathfrak u}}
\newcommand{\tuu}{\tilde{\mathfrak u}}
\newcommand{\buu}{\bar{\mathfrak u}}
\newcommand{\bnn}{\bar{\mathfrak n}}
\newcommand{\bab}{\bar{\mathfrak b}}
\newcommand{\bns}{\bar{\mathfrak n}_{\T}}
\newcommand{\mmm}{{\mathfrak m}}
\newcommand{\smm}{\mbox{}^s{\mathfrak m}}
\newcommand{\mms}{{\mathfrak m}_{\T}}
\newcommand{\jjj}{{\mathfrak j}}
\newcommand{\hhh}{{\mathfrak h}}
\newcommand{\shh}{\mbox{}^s{\mathfrak h}}
\newcommand{\uhh}{\mbox{}^u{\mathfrak h}}
\newcommand{\hhd}{{\mathfrak h}^\ast}
\newcommand{\qqq}{{\mathfrak q}}
\newcommand{\ppp}{{\mathfrak p}}
\newcommand{\tpp}{\tilde{\mathfrak p}}
\newcommand{\tnn}{\tilde{\mathfrak n}}
\newcommand{\tqq}{\tilde{\mathfrak q}}
\newcommand{\pps}{{\mathfrak p}_{\T}}
\newcommand{\bapp}{\bar{\mathfrak p}}
\newcommand{\baps}{\bar{\mathfrak p}_{\T}}
\newcommand{\sss}{{\mathfrak s}}
\newcommand{\ooo}{{\mathfrak o}}
\newcommand{\ddd}{{\mathfrak d}}
\newcommand{\Sth}{{\mathfrak S}}

\newcommand{\ggc}{{\mathfrak g}_{\sf c}}
\newcommand{\bbc}{{\mathfrak b}_{\sf c}}
\newcommand{\llc}{{\mathfrak l}_{\sf c}}
\newcommand{\uuc}{{\mathfrak u}_{\sf c}}
\newcommand{\nsc}{{({\mathfrak n}_{\T})}_{\sf c}}
\newcommand{\nnc}{{({\mathfrak n}_{\T})}_{\sf c}}
\newcommand{\buc}{\bar{\mathfrak u}_{\sf c}}
\newcommand{\bnc}{{(\bar{\mathfrak n}_{\T})}_{\sf c}}
\newcommand{\hhc}{{\mathfrak h}_{\sf c}}
\newcommand{\hcd}{{{\mathfrak h}_{\sf c}}^\ast}
\newcommand{\kkc}{{\mathfrak k}_{\sf c}}
\newcommand{\ppc}{{\mathfrak p}_{\sf c}}
\newcommand{\pcs}{({\mathfrak p}_{\T})_{\sf c}}
\newcommand{\ncs}{({\mathfrak n}_{\T})_{\sf c}}
\newcommand{\bpp}{\bar{\mathfrak p}}
\newcommand{\bqq}{\bar{\mathfrak q}}

\newcommand{\spl}{\mbox{}^s}

\newcommand{\gl}{{\mathfrak g}{\mathfrak l}}
\newcommand{\gll}{{\mathfrak gl}_L}
\newcommand{\glr}{{\mathfrak gl}_R}
\newcommand{\so}{\sss\ooo}

\newcommand{\sA}{{\mbox{}^s{A}}}
\newcommand{\sM}{{\mbox{}^s{M}}}
\newcommand{\sT}{{\mbox{}^s{T}}}
\newcommand{\sH}{{\mbox{}^s{H}}}
\newcommand{\sB}{{\mbox{}^s{P}}}
\newcommand{\sP}{{\mbox{}^s{P}}}
\newcommand{\uH}{{\mbox{}^u{H}}}

\newcommand{\bi}{{\bold i}}
\newcommand{\fsp}{{\mathfrak s}{\mathfrak p}}

\newcommand{\Gc}{{G_{\mathbb C}}}
\newcommand{\Gac}{{G_{\mathbb C}^{\rm ad}}}
\newcommand{\Gca}{{G_{\mathbb C}^{\rm ad}}}
\newcommand{\Gf}{{G^{\flat}}}
\newcommand{\Gs}{{G^{\sharp}}}
\newcommand{\Ga}{{G^{\ast}}}
\newcommand{\Kc}{{K_{\mathbb C}}}
\newcommand{\Kcf}{{K_{\mathbb C}^\flat}}
\newcommand{\Kf}{{K^\flat}}
\newcommand{\Ks}{{K^\sharp}}
\newcommand{\Ka}{{K^\ast}}
\newcommand{\Kcs}{{K_{\mathbb C}^\sharp}}
\newcommand{\Kca}{{K_{\mathbb C}^\ast}}
\newcommand{\Bc}{{B_{\mathbb C}}}
\newcommand{\bBc}{{\bar{B}_{\sf C}}}
\newcommand{\Pc}{{P_{\mathbb C}}}
\newcommand{\bPc}{{\bar{P}_{\mathbb C}}}
\newcommand{\Nc}{{N_{\mathbb C}}}
\newcommand{\bNc}{{\bar{N}_{\mathbb C}}}

\newcommand{\hol}{{\cal O}}
\newcommand{\dif}{{\cal D}}
\newcommand{\ana}{{\cal A}}

\newcommand{\tpca}{\tilde{\cal P}}
\newcommand{\aca}{{\cal A}}
\newcommand{\bca}{{\cal B}}
\newcommand{\cca}{{\cal C}}
\newcommand{\dca}{{\cal D}}
\newcommand{\diff}{{\cal D}}
\newcommand{\eca}{{\cal E}}
\newcommand{\fca}{{\cal F}}
\newcommand{\gca}{{\cal G}}
\newcommand{\hca}{{\cal H}}
\newcommand{\ica}{{\cal I}}
\newcommand{\jca}{{\cal J}}
\newcommand{\kca}{{\cal K}}
\newcommand{\lca}{{\cal L}}
\newcommand{\mca}{{\cal M}}
\newcommand{\nca}{{\cal N}}
\newcommand{\oca}{{\cal O}}
\newcommand{\pca}{{\cal P}}
\newcommand{\qca}{{\cal Q}}
\newcommand{\rca}{{\cal R}}
\newcommand{\sca}{{\cal S}}
\newcommand{\tca}{{\cal T}}
\newcommand{\uca}{{\cal U}}
\newcommand{\vca}{{\cal V}}
\newcommand{\wca}{{\cal W}}
\newcommand{\xca}{{\cal X}}
\newcommand{\yca}{{\cal Y}}
\newcommand{\zca}{{\cal Z}}

\newcommand{\HG}{{\cal H}_G}
\newcommand{\PR}{{\cal P}r}
\newcommand{\gS}{{\mathfrak S}}
\newcommand{\ii}{\sqrt{-1}}
\newcommand{\real}{\mbox{Re}}
\newcommand{\Real}{\mbox{Re}}
\newcommand{\res}{\mbox{res}}
\newcommand{\supp}{\mbox{supp}}
\newcommand{\rank}{\mbox{rank}}
\newcommand{\card}{\mbox{card}}
\newcommand{\Ad}{\mbox{Ad}}
\newcommand{\ad}{\mbox{ad}}
\newcommand{\dg}{{\deg}}
\newcommand{\Hom}{\mbox{Hom}}
\newcommand{\End}{\mbox{End}}
\newcommand{\Ext}{\mbox{Ext}}
\newcommand{\fgt}{\mbox{Fgt}}
\newcommand{\JH}{\mbox{JH}}
\newcommand{\Dim}{\mbox{Dim}}
\newcommand{\gr}{\mbox{gr}}
\newcommand{\Ass}{\mbox{Ass}}
\newcommand{\WF}{\mbox{WF}}
\newcommand{\AS}{\mbox{AS}}
\newcommand{\Ann}{\mbox{Ann}}
\newcommand{\LAnn}{\mbox{LAnn}}
\newcommand{\RAnn}{\mbox{RAnn}}
\newcommand{\tr}{\mbox{tr}}
\newcommand{\hht}{\mbox{ht}}
\newcommand{\Mod}{{\sf Mod}}
\newcommand{\sgn}{\mbox{sgn}}
\newcommand{\pro}{\mbox{pro}}
\newcommand{\Ind}{\mbox{\sf Ind}}
\newcommand{\SO}{\mbox{SO}}
\newcommand{\Oo}{\mbox{O}}
\newcommand{\SOo}{{\mbox{SO}_0}}
\newcommand{\GL}{\mbox{GL}}
\newcommand{\Spp}{\mbox{Sp}}
\newcommand{\U}{\mbox{U}}
\newcommand{\triv}{{\mbox{triv}}}
\newcommand{\Rea}{{\mbox{Re}}}

\newcommand{\Res}{\mbox{Res}}
\newcommand{\diag}{\mbox{diag}}

\newcommand{\PP}{{\sf P}}
\newcommand{\PS}{{\sf P}_{\T}^{++}}
\newcommand{\PSP}{{\sf P}_{\T^\prime}^{++}}
\newcommand{\PPS}{{\sf P}_{\T}^{++}}
\newcommand{\PPSA}{{\sf P}_{\T^\alpha}^{++}}
\newcommand{\PPSK}{{\sf P}_{\T^k}^{++}}

\newcommand{\Q}{{\sf Q}}

\newcommand{\ww}{{\sf w}_\psi}
\newcommand{\Wh}{{\sf Wh}_{\bnn,\psi}^\infty}
\newcommand{\Whc}{{\sf Wh}_{\bnn,\psi}^\infty}
\newcommand{\ind}{{\rm ind}}
\newcommand{\hind}{\mbox{\rm h-ind}}
\newcommand{\Whd}{{\sf Wh}_\Psi^\ast}
\newcommand{\HC}{Harish-Chandra $(\ggc,\kkc)$-module}
\newcommand{\HCs}{Harish-Chandra $(\ggc,\kkc)$-modules}
\newcommand{\cH}{\check{H}}
\newcommand{\bh}{\bar{h}}

\newcommand{\lel}{\stackrel{L}{\leq}}
\newcommand{\lell}{\stackrel{L}{\leq}}
\newcommand{\ler}{\stackrel{R}{\leq}}
\newcommand{\lerr}{\stackrel{R}{\leq}}
\newcommand{\lelr}{\stackrel{LR}{\leq}}

\newcommand{\gel}{\stackrel{L}{\geq}}
\newcommand{\ger}{\stackrel{R}{\leq}}
\newcommand{\gelr}{\stackrel{LR}{\leq}}

\newcommand{\lleq}{\stackrel{L}{\sim}}
\newcommand{\rreq}{\stackrel{R}{\sim}}
\newcommand{\lreq}{\stackrel{LR}{\sim}}

\newcommand{\leqs}{\leqslant}
\newcommand{\geqs}{\geqslant}
\title{{\sf On the homomorphisms between scalar generalized Verma modules}}
\author{{\sf  Hisayosi Matumoto}
\\Graduate School of Mathematical Sciences\\ University of Tokyo\\ 3-8-1
Komaba, Tokyo\\ 153-8914, JAPAN\\ e-mail: hisayosi@ms.u-tokyo.ac.jp}
\date{}
\maketitle
\begin{abstract}
In this article, we study   the homomorphisms between scalar generalized
Verma modules.  We conjecture that any homomorphism between  is composition of elementary homomorphisms.  The purpose of this article is to 
 show the conjecture is affirmative for  many  parabolic subalgebras under the assumption that the infinitesimal characters are regular.
\footnote{This work was supported by Grant-in-Aid for Scientific Research (No.\ 2054001)}
\end{abstract}

\setcounter{section}{0}
\section*{\S\,\, 0.\,\,\,\, Introduction}
\setcounter{subsection}{0}

We study   the homomorphisms between generalized
Verma modules, which are induced from one dimensional representations
(such generalized Verma modules are called scalar, cf.\   \cite{[Boe]}).

 Classification of the homomorphisms between scalar generalized Verma modules  is equivalent to that of equivariant differential
operators between the spaces of sections of homogeneous line bundles on
generalized flag manifolds. (cf.\ \cite{[Ko]},\cite{[Do]}, \cite{ [Jk]}, \cite{
[CS]}, and \cite{[Hu]}.)

In \cite{[Vr]}, Verma constructed  homomorphisms between
Verma modules associated with root reflections.
Bernstein, I.\ M.\ Gelfand, and S.\ I.\ Gelfand proved that all the nontrivial homomorphisms between Verma modules are compositions of homomorphisms constructed by Verma.  (\cite{[BGG]}) 

Later, Lepowsky studied the generalized Verma modules.
In particular, Lepowsky (\cite{[L2]}) constructed a class
of nontrivial homomorphisms between scalar generalized Verma modules
associated to the parabolic
subalgebras which are the complexifications of the minimal parabolic
subalgebras of real reductive Lie algebras. 

In \cite{[M]}, elementary homomorphisms between scalar generalized Verma modules
are introduced.  They can be regarded as a generalization of homomorphisms introduced by Verma and Lepowsky.

We propose a conjecture on the classification of the homomorphisms between scalar generalized Verma modules, which can be regarded as a generalization of the above-mentioned result of Bernstein-Gelfand-Gelfand.　　

{\bf Conjecture A}  \,\,\,\, All the nontrivial homomorphisms between scalar generalized Verma modules are compositions of elementary homomorphisms.
\medskip

Soergel's result (\cite{[So]} Theorem 11) implies that Conjecture A is reduced to the integral infinitesimal character case.

The purpose of this article is to show the conjecture is affirmative for  many  parabolic subalgebras under the assumption that the infinitesimal characters are regular.
In order to explain our results, we introduce some notations.
Let $\ggg$ be a complex reductive Lie algebra and we fix  a Cartan subalgebra $\hhh$ of $\ggg$.  We denote by $\Delta$ (resp.\ $W$) the root system (resp.\ the Weyl group) with respect to
$(\gggg,\hhh)$.
We fix a basis $\Pi$ of $\Delta$.
For $\T\subsetneq\Pi$, we put $\aas=\{H\in\hhh\mid \forall \alpha\in \T \,\,\,\alpha(H)=0\}$ and $\Sigma_\T=\{\alpha|_{\aas}\mid\alpha\in\Delta\}-\{0\}$.
We denote by $\ppp_\T$ the standard parabolic subalgebra corresponding to $\T$ and by $\llll_\T$ its Levi subalgebra containing $\hhh$.  We consider the Weyl group for parabolic sub algebra $W(\T)=\{w\in\mid w\T=\T\}$.

We call $\T$ normal if any two parabolic subalgebras with the Levi part $\llll_\T$ are conjugate under an inner automorphism of $\ggg$. If $\T$ is normal, we call $\ppp_\T$ normal.
For example, a complexified minimal parabolic subalgebras of real simple Lie algebras except $\sss\uuu(p,q)$ \,\,\, $(p-1>q>0)$, $\sss\ooo^\ast(4n+2)$, $\eee_{6(-14)}$ are normal.
Roughly speaking, if $\T$ is normal, the reflection $\sigma_\gamma$ on $\aas$ with respect to $\gamma\in \Sigma_\T$ can be regarded as an involution of the Weyl group for $(\ggg,\hhh)$.
A normal subset $\T$ of $\Pi$ is called strictly normal, if $\sigma_\gamma$ is a Duflo involution of some Weyl group (see Definition 4.2.1) .
If $\T$ is strictly normal, there exists an elementary homomorphism with repect to $\sigma_\gamma$ for each $\gamma\in\Sigma_\T$.

Let $\ppp_\T$ be  a complexified minimal parabolic subalgebra of a  real simple Lie algebra and assume $\ppp_\T$ is normal but is not strictly normal.    Then, $\ppp_\T$ is a complexified minimal parabolic subalgebra of $\sss\ooo(2n+1-q,q)$ \,\,\, $(n>q\geqslant 1)$, or $\sss\ppp(n,n)$  \,\,\, $(n\geqslant 1)$.

One of the main result of this article is the following theorem.

{\bf Theorem B} (Theorem 5.1.3) \,\,\,\, {\it If $\T$ is strictly normal, then each nontrivial homomorphism between scalar generalized Verma modules induced from  $\ppp_\T$ with regular integral  infinitesimal character is composition of elementary homomorphisms.}
\medskip

We also have the following result. 

{\bf Theorem C} (cf.\ Theorem 5.1.3, Corollary 6.8.7, Proposition 6.5.2) \,\,\,\,  {\it If $\T$ is normal and $\ggg$ is an exceptional Lie algebra, then each nontrivial homomorphism between scalar generalized Verma modules induced from  $\ppp_\T$ with a regular integral  infinitesimal character is composition of elementary homomorphisms.}
\medskip

For $\ggg\llll(n,\cpx)$, we have the following result.  (In fact, we have a more general result.)

{\bf Theorem D} (Corollary 7.3.6) \,\,\,\,  {\it If $\ppp_\T$ is a complexified minimal parabolic subalgebra of $\sss\uuu(p,q)$, then each nontrivial homomorphism between scalar generalized Verma modules induced from  $\ppp_\T$ with a regular integral  infinitesimal character is composition of elementary homomorphisms.}
\medskip

Finally,  we have the following result.

{\bf Corollary  E} (Theorem 5.1.3, Corollary 6.3.7) \,\,\,\, {\it 
Let $\ppp_\T$ be  a complexified minimal parabolic subalgebra of a real simple Lie algebra other than   $\sss\ooo^\ast(4n+2)$,  $\eee_{6(-14)}$, $\sss\ooo(2n+1-q,q)$ \,\,\, $(n>q>2)$, and $\sss\ppp(n,n)$  \,\,\, $(n>1)$.
Then, each nontrivial homomorphism between scalar generalized Verma modules induced from  $\ppp_\T$ with regular integral  infinitesimal character is composition of elementary homomorphisms.} 

For $\sss\ooo^\ast(4n+2)$ and  $\eee_{6(-14)}$, we show a weaker statement holds.  (cf.\ 8.3, 8.4)

This article consists of nine sections.

We fix notations and introduce some fundamental material in \S 1.

In \S 2, we explain how to reduce the problem to the integral infinitesimal character case.  We also show that we can associate an element of  $W(\T)$ to a homomorphism between generalized Verma modules with regular infinitesimal characters.  Finally we formulate the translation principle for generalized Verma modules, which is essentially obtained in \cite{[VU]}, \cite{[VI]}.

In \S 3, we introduce the notion of normal parabolic subalgebras and describe the classification.  We  prove that the Bruhat ordering on $W(\T)$ coincides with the restriction of that of $W$ to $W(\T)$ for each normal $\T$.

In \S 4, we introduce the notion of an elementary homomorphism and describe related notions and results.

In \S 5, we introduce the notion of strictly normal parabolic subalgebras and describe the classification.  We also prove Theorem B.

In \S 6, we consider normal but not strictly normal case.  We prove that our conjecture is affirmative for a complexified minimal parabolic subalgebras of $\sss\ooo(2n-1,2)$, and  normal but strictly normal parabolic subalgebras of exceptional algebras in regular integral infinitesimal character case.  

In \S 7, We consider $\ggg\llll(n,\cpx)$.  We introduce the notion of almost normal parabolic subalgebras.  We prove that our conjecture is affirmative for an almost normal parabolic subalgebras of $\ggg\llll(n,\cpx)$ in regular integral infinitesimal character case.  

In \S 8, we consider complexified minimal parabolic subalgebras of $\sss\ooo^\ast(4n+2)$ and $\eee_{6(-14)}$.

In \S 9, we treat two examples.

\setcounter{section}{1}
\setcounter{subsection}{0}

\section*{\S\,\, 1.\,\,\,\,Notations and Preliminaries }

\subsection{General notations}

In this article, we use the following notations and conventions.

As usual we denote the complex number field, the real number field, the
ring of (rational) integers, and the set of non-negative integers by
$\cpx$, $\rel$, $\itg$, and $\nat$ respectively.
$\frac{1}{2}\nat$ means the set $\left\{\left. \frac{n}{2}\right|
n\in\nat\right\}$, and
$\frac{1}{2}+\nat$ means the set $\left\{\left. \frac{1}{2}+n\right| n\in\nat\right\}$.
We denote by $\emptyset$ the empty set.
For any (non-commutative) $\cpx$-algebra $R$, ``ideal'' means  ``2-sided ideal'', ``$R$-module'' means ``left $R$-module'', and sometimes we denote by $0$ (resp.\ $1$) the trivial $R$-module $\{0\}$ (resp.\ $\cpx$).
Often, we identify a (small) category and the set of its objects.
Hereafter ``$\dim$'' means the dimension as a complex vector space, and ``$\otimes$'' (resp. $\Hom$) means the tensor product over $\cpx$ (resp. the space of $\cpx$-linear mappings), unless we specify.
For a complex vector space $V$, we denote by $V^\ast$ the dual vector space.
For $a,b\in\cpx$, ``$a\leqslant b$'' means that $a,b\in\rel$ and $a\leqslant b$.
We denote by $A- B$ the set theoretical difference.
$\card A$ means the cardinality of a set $A$.

\subsection{Notations for reductive Lie algebras}

Let $\gggg$ be a complex reductive Lie algebra, $U(\gggg)$ the universal
enveloping algebra of $\gggg$, and  $\hhh$ a Cartan subalgebra of
$\gggg$.
We denote by $\Delta$ the root system with respect to
$(\gggg,\hhh)$.
We fix some positive root system $\Delta^+$ and let $\Pi$ be the set of simple roots.
Let $W$ be the Weyl group of the pair $(\gggg, \hhh)$ and let
$\langle\,\,,\,\,\rangle$ be a non-degenerate invariant bilinear form on
$\gggg$.
For $w\in W$, we denote by $\ell(w)$ the length of $w$ as usual.
We also denote the inner product on $\hhd$ which is induced from the
above form by the same symbols $\langle\,\,,\,\,\rangle$.
For $\alpha\in\Delta$, we denote by $s_\alpha$ the reflection in $W$ with respect to $\alpha$.
We denote by $w_0$ the longest element of $W$.
For $\alpha\in\Delta$, we define the coroot ${\alpha}^\vee$ by
${\alpha}^\vee = \frac{2\alpha}{\langle\alpha,\alpha\rangle}$,
as usual.  We denote by $\Delta^\vee$ the dual root system $\{\alpha^\vee\mid \alpha\in\Delta\}$.
We call $\lambda\in \hhd$ is dominant (resp. anti-dominant), if $\langle\lambda,{\alpha}^\vee\rangle$ is not a negative (resp. positive) integer, for each $\alpha\in \Delta^+$.
We call $\lambda\in\hhd$ regular, if $\langle\lambda,\alpha\rangle\neq 0$, for each $\alpha\in\Delta$.
We denote by $\PP$ the integral weight lattice, namely
$\PP = \{\lambda\in\hhd\mid\langle\lambda,{\alpha}^\vee\rangle\in\itg \,\,\,\,\hbox{for all}\,\,\, \alpha \in\Delta\}$.
If $\lambda\in\hhd$ is contained in $\PP$, we call $\lambda$ an integral weight.
We define $\rho\in\PP$ by $\rho = \frac{1}{2}\sum_{\alpha\in\Delta^+}\alpha$.
Put 
$\gggg_\alpha =\{X\in\gggg\mid\forall H\in\hhh \;\; [H,X] =\alpha(H)X\}$,
$\uuu = \sum_{\alpha\in\Delta^+}\gggg_\alpha$,
$\bbb = \hhh+\uuu$.
Then $\bbb$ is a Borel subalgebra of $\gggg$.
We denote by $\Q$ the root lattice, namely $\itg$-linear span of $\Delta$.  We also denote by $\Q^+$ the linear combination of $\Pi$ with non-negative integral coefficients.  
 For $\lambda\in\hhh^\ast$, we denote by $W_\lambda$ the integral Weyl group.
Namely,
\[ W_\lambda=\{w\in W\mid w\lambda-\lambda\in\Q\}.\]
We denote by $\Delta_\lambda$ the set of integral roots.
\[\Delta_\lambda=\{\alpha\in\Delta\mid \langle\lambda,{\alpha}^\vee\rangle\in\itg\}.\]
It is well-known that $W_\lambda$ is the Weyl group for $\Delta_\lambda$.
We put $\Delta^+_\lambda=\Delta^+\cap\Delta_\lambda$.
This is a positive system of $\Delta_\lambda$.  We denote by $\Pi_\lambda$ the set of simple roots for $\Delta^+_\lambda$ and denote by $S_\lambda$  (resp.\ $S$) the set of reflection corresponding to the elements in $\Pi_\lambda$ (resp.\ $\Pi$).  
So, $(W_\lambda, S_\lambda)$  and $(W, S)$ are Coxeter systems.
We denote by $\Q_\lambda$ the integral root lattice, namely
$\Q_\lambda=\itg\Delta^+_\lambda$ and put
$\Q^+_\lambda=\nat\Pi_\lambda$.

Next, we fix notations for a parabolic subalgebra (which contains $\bbb$).
Hereafter, through this article we fix an arbitrary subset $\T$ of $\Pi$.
Let $\langle\T\rangle$ be the set of the elements of $\Delta$ which are written by linear combinations of elements of $\T$ over $\itg$.
Put
$\aas = \{H\in\hhh\mid \forall \alpha\in \T \,\,\,\alpha(H)=0\}$,
$\lls = \hhh +\sum_{\alpha\in\langle{\T}\rangle}\gggg_\alpha$,  
$\nns = \sum_{\alpha\in\Delta^+-\langle{\T}\rangle}\gggg_\alpha$,  
$\pps = \lls +\nns$.
Then $\pps$ is a parabolic subalgebra of $\gggg$ which contains $\bbb$.
Conversely, for an arbitrary parabolic subalgebra $\ppp\supseteq\bbb$, there exists some $\T \subseteq \Pi$ such that $\ppp =\pps$.
We denote by $W_\T$ the Weyl group for $(\lls,\hhh)$. 
$W_\T$ is identified with a subgroup of $W$ generated by $\{s_\alpha\mid \alpha\in\T\}$.
We denote by $w_\T$ the longest element of $W_\T$.
Using the invariant non-degenerate bilinear form $\langle\,\,,\,\,\rangle$, we regard  ${\aaa_{\T}}^\ast$ as a subspace of $\hhh^\ast$. 

Put
$\rho_{\T} = \frac{1}{2}(\rho-w_{\T}\rho)$ and 
$\rho^{\T} = \frac{1}{2}(\rho+w_{\T}\rho)$.
Then, $\rho^{\T}\in{\aaa_{\T}}^\ast$.

For $\T\subsetneq\Pi$, we define "the restricted root system" as follows.
\[ \Sigma_\T=\{\alpha|_{\aaa_\T}\mid\alpha\in\Delta\}-\{0\}.\]
\[ \Sigma_\T^+=\{\alpha|_{\aaa_\T}\mid\alpha\in\Delta^+\}-\{0\}.\]
Unfortunately, in general, $ \Sigma_\T$ does not satisfy the axioms of the root systems.

Define
\begin{align*}
\PS & = \{\lambda\in\hhd\mid \forall\alpha\in \T\,\,\,\,\,
 \langle\lambda,{\alpha}^\vee\rangle\in\{1,2,...\}\}\\
{}^\circ\PS & =\{\lambda\in\hhd\mid \forall\alpha\in \T\,\,\,\,\,
 \langle\lambda,{\alpha}^\vee\rangle=1\}
\end{align*}
We easily have
\begin{align*}
{}^\circ\PS =\{\rho_{\T}+\mu\mid \mu\in\ads\}.
\end{align*}
For $\mu\in\hhd$ such that $\mu+\rho\in\PS$, we denote by $\sigma_\T(\mu)$ the irreducible finite-dimensional $\lls$-representation whose highest weight is $\mu$.
Let $E_\T(\mu)$ be the representation space of $\sigma_\T(\mu)$.
We define a left action of $\nns$ on $E_\T(\mu)$ by $X\cdot v =0$ for all $X\in\nns$ and $v\in E_\T(\mu)$.
So, we regard $E_\T(\mu)$ as a $U(\pps)$-module.

For $\mu\in\PS$, we define a generalized Verma module (\cite{[L4]}) as follows.
\[ M_\T(\mu) = U(\gggg)\otimes_{U(\pps)}E_\T({\mu-\rho}).\]
For all $\lambda\in\hhd$, we write $M(\lambda) = M_\emptyset(\lambda)$.
$M(\lambda)$ is called a Verma module.
For $\mu\in\PS$, $M_\T(\mu)$ is a quotient  module of $M(\mu)$.
Let $L(\mu)$ be the unique highest weight $U(\gggg)$-module with the highest weight $\mu-\rho$.
Namely, $L(\mu)$ is a unique irreducible quotient of $M(\mu)$.
For $\mu\in\PS$, the canonical projection of $M(\mu)$ to
$L(\mu)$ is factored by $M_\T(\mu)$.

$\dim E_{\T}(\mu-\rho)=1$ if and only if $\mu\in {}^\circ\PS$.
If $\mu\in {}^\circ\PS$, we call $M_\T(\mu)$ a scalar generalized
Verma module.

\subsection{Translation functors}

We denote by $Z(\ggg)$ the center of $U(\ggg)$.
It is well-known that $Z(\ggg)$ acts on $M(\lambda)$ by the Harish-Chandra homomorphism $\chi_\lambda : Z(\ggg) \rightarrow \cpx$ for all $\lambda$.
$\chi_\lambda =\chi_\mu$ if and only if there exists some $w\in W$ such that $\lambda =w\mu$.
We denote by $\hbox{\bf Z}_\lambda$ the kernel of $\chi_\lambda$ in $Z(\ggg)$.
Let $M$ be a $U(\ggg)$-module and $\lambda\in\hhh^\ast$.
We say that $M$ has an infinitesimal character $\lambda$ iff $Z(\ggg)$ acts on $M$ by $\chi_\lambda$.
We say that $M$ has a generalized infinitesimal character $\lambda$ iff for any $v\in M$ there is some positive integer $n$ such that ${\hbox{\bf Z}_\lambda}^nv=0$.
  We say $M$ is locally $Z(\ggg)$-finite, iff for any $v\in M$ we have  $\dim Z(\ggg)v<\infty$.
We denote by $\mca_{Zf}$ (cf \cite{[BG]}) the category of $Z(\ggg)$-finite $U(\ggg)$-modules.
We also denote by $\mca[\lambda]$ the category of $U(\ggg)$-modules with generalized infinitesimal character $\lambda$.
Then, from the Chinese remainder theorem, we have a direct sum of abelian categories 
$\mca_{Zf}=\bigoplus_{\lambda\in\hhh^\ast}\mca[\lambda]$.
We denote by $P_\lambda$ the projection functor from $\mca_{Zf}$ to $\mca[\lambda]$.
For $\mu\in\PP$, we denote by $V_\mu$ the irreducible finite-dimensional $U(\ggg)$-module with an extreme weight $\mu$.
Let $\mu,\lambda\in\hhh^\ast$ satisfy $\mu-\lambda\in\PP$.  
Let $M$ be an object of $\mca[\lambda]$.
Then, from a result of Kostant we have that $M\otimes V_{\mu-\lambda}$ is an object of $\mca_{Zf}$.  
So, we can define a translation functor $T^\mu_\lambda$ from $\mca[\lambda]$ to $\mca[\mu]$ as follows.
\[T^\mu_\lambda(M)=P_\mu(M\otimes V_{\mu-\lambda}).\]
$T^\mu_\lambda$ is an exact functor.

\subsection{Submodules of scalar generalized Verma modules}
For a finitely generated $U(\ggg)$-module $V$, we denote by $\Dim(V)$ the Gelfand-Kirillov dimension of $V$  (cf.\ \cite{[Vg]}).

The following proposition is more or less known.
\begin{prop}
Let $\T\subseteq\Pi$.
Then we have :
\begin{itemize}
\item[(1)] Let $\lambda\in\PS$.  Then, for each nonzero submodule $X$ of $M_\T(\lambda)$, we have $\Dim(X)=\Dim( M_\T(\lambda))=\dim\nns$.
\item[(2)] Let $\lambda\in{}^\circ\PS$.   Then, for each nonzero submodule $X$ of $M_\T(\lambda)$, we have $\Dim(M_\T(\lambda)/X)<\Dim( M_\T(\lambda))=\dim\nns$.
\item[(3)]  Let $\lambda\in{}^\circ\PS$. Then $M_\T(\lambda)$ has a unique irreducible submodule.
\end{itemize}
\end{prop}
\proof
Since $M_\T(\lambda)$ is free of finite rank as a $U(\bns)$-module, we have 
$\Dim( M_\T(\lambda))=\dim\bns=\dim\nns$.
A nonzero submodule $X$ of $M_\T(\lambda)$ is torsion free as $U(\bns)$-module, so $\dim\bns\leqslant\Dim(X)\leqslant \Dim( M_\T(\lambda))=\dim\bns$.
So, we have (1).
Next, let $\lambda\in{}^\circ\PS$.
Then, the multiplicity (the Bernstein degree) (cf.\ \cite{[Vg]}) of $ M_\T(\lambda)$ is one.  So, the number of the irreducible irreducible constituents of $M_\T(\lambda)$ which have the Gelfand -Kirillov dimension $\dim\nns$ is one.
So, from (1), we have (2) and (3).  \,\,\,\, $\Box$

\setcounter{section}{2}
\setcounter{subsection}{0}

\section*{\S\,\, 2.\,\,\,\, Formulation of the problem }

We retain the notation of \S 1.
In particular, $\T$ is a proper subset of $\Pi$.

\subsection{Basic results of Lepowsky}

The following result is one of the fundamental results on the existence
problem of homomorphisms between scalar generalized Verma modules.

\begin{thm} (\cite{[L3]})

Let $\mu,\nu\in {}^\circ\PS$.

(1) \,\,\,\, $\dim \Hom_{U(\gggg)}(M_{\T}(\mu), M_{\T}(\nu))\leqslant 1.$

(2) \,\,\,\, Any non-zero homomorphism of $M_{\T}(\mu)$ to
 $M_{\T}(\nu)$ is injective.
\end{thm}

Hence, the classification problem of homomorphisms between generalized Verma modules is reduce to the following problem.

{\bf Problem 1} \,\, Let $\mu,\nu\in {}^\circ\PS$.
 When is $M_{\T}(\mu)\subseteq M_{\T}(\nu)$ ?

\subsection{Reduction to the integral infinitesimal character setting}

Since the both $\nu\in W\mu$ and  $\nu-\mu\in Q^+$ are necessary
conditions for $M_{\T}(\mu)\subseteq M_{\T}(\nu)$, we can reformulate our problem as
follows.

{\bf Problem 2} \,\, Let $\lambda\in {}^\circ\PS$ be dominant.
Let $x,y\in W_\lambda$ be such that $x\lambda,y\lambda\in{}^\circ\PS$.
 When is $M_{\T}(x\lambda)\subseteq M_{\T}(y\lambda)$ ? 

We fix $\lambda\in {}^\circ\PS$.
Then,  we can construct a suralgebra $\gggg^\prime$
 of $\hhh$ such that the corresponding Coxeter system is $(W_\lambda,\Phi_\lambda)$.
Since $\T\subsetneq\Pi_\lambda$ holds, we can construct the corresponding
 parabolic subalgebra $\ppp_{\T}^\prime$ of $\gggg^\prime$.
For $\mu\in\PS$, we denote by $M^\prime_{\T}(\mu)$ the corresponding
generalized Verma module of $\gggg^\prime$.
We consider the category $\hol$ in the sense of
\cite{[BGG]} corresponding to our particular choice
of positive root system.
More precisely, we denote by $\hol$ (respectively $\hol^\prime$) ``the
category $\hol$'' for $\gggg$ (respectively $\gggg^\prime$).
We denote by $\hol_\lambda$ (respectively,  $\hol_\lambda^\prime$) the
full subcategory of $\hol$ (respectively $\hol^\prime$) consisting of
the objects with a generalized infinitesimal character $\lambda$.
Soergel's celebrated theorem (\cite{[So]} Theorem 11) says that there is
a Category equivalence between $\hol_\lambda$ and $\hol^\prime_\lambda$.
Under the equivalence a Verma module $M(x\lambda)$ \,\,\, $(x\in
W_\lambda)$ corresponds to $M^\prime(x\lambda)$.
From Lepowsky's generalized BGG resolutions of the generalized Verma
modules and their rigidity, we easily see $M_{\T}(x\lambda)$
corresponds to $M^\prime_{\T}(x\lambda)$ under Soegel's category
equivalence.
So, we have the following lemma as a corollary of Soergel's theorem.

\begin{lem}  
Let $\lambda\in \hhh^\ast$ be dominant.
Let $x,y\in W_\lambda$ be such that $x\lambda,y\lambda\in{}^\circ\PS$.
Then, the following two conditions are equivalent.

(1) \,\, $M_{\T}(x\lambda)\subseteq M_{\T}(y\lambda)$.

(2) \,\, $M_{\T}^\prime(x\lambda)\subseteq M_{\T}^\prime(y\lambda)$.
\end{lem}

This lemma tells us that  we may reduce Problem 2 to the case that  $\lambda$
is integral. 

We discuss  another application of Soergel's theorem.
We denote by $\ggg^\vee$ the reductive Lie algebra corresponding to the coroot system $\Delta^\vee$.  We regard a Cartan subalgebra $\hhh$ as a Cartan subalgebra of $\ggg^\vee$.
We attach $\vee$ to the notion with respect to $\ggg^\vee$ corresponding to that of $\ggg$.  Then we have canonical isomorphism $(W, S)\cong (W^\vee,S^\vee)$ of the Coxeter systems.  So, we identify them.
For $\T\subsetneq\Pi$, we put $\T^\vee=\{\alpha^\vee\mid\alpha\in\T\}\subsetneq\Pi^\vee$.
We put ${}^\circ\PP_{\T^\vee}^{\vee++}=\{\lambda\in\hhd\mid \langle\lambda,\alpha\rangle=1 \,\,\, (\alpha\in\T)\}$.
For $\lambda\in {}^\circ\PP_{\T^\vee}^{\vee++}$, we consider a scalar generalized Verma module $M_{\T^\vee}^\vee(\lambda)$ of $\ggg^\vee$.
The following result is an immediate consequence of Soergel's theorem.
\begin{thm}
Let $\lambda\in\PP$ and  $\mu\in\PP^\vee$ be dominant regular.
Let $x,y\in W=W^\vee$.  We assume that $x\lambda,y\lambda\in {}^\circ\PS$ and $x\mu,y\mu\in {}^\circ\PP_{\T^\vee}^{\vee++}$.
Then, $M_\T(x\lambda)\subseteq M_\T(y\lambda)$ if and only if $M_{\T^\vee}^\vee(x\mu)\subseteq M_{\T^\vee}^\vee(y\mu)$.
\end{thm}
Hence, we may reduce Problem 1 for simple Lie algebras of the type $\hbox{\rm C}_n$ to that for simple Lie algebras of the type $\hbox{\rm B}_n$.

\subsection{Comparison of $\tau$-invariants}

We put 
\[ W({\T})=\{w\in W\mid w\T= \T \}.\]

Then, we easily  have the following lemma.

\begin{lem}
We have

(a) \,\,\, $W({\T})=\{w\in W\mid \mbox{$w{\alpha}^\vee\in\T^\vee$ for all $\alpha\in\T$.} \}$.

(b) \,\,\, 
$W({\T})=\{w\in W\mid w\rho_{\T}=\rho_{\T}, w\T\subseteq\Delta^+\}$.
 
(c) \,\,\,  $w_{\T}w=ww_{\T}$ for all $w\in W({\T})$.

(d) \,\,\, $W({\T})$ preserves $\aas^\ast$.

(e) \,\,\, $W({\T})\subseteq W_{\rho_{\T}}$.
\end{lem}

In this section, we prove the following proposition.

\begin{prop}
Let $\lambda\in {}^\circ\PS$ be regular.
Let $x\in W_\lambda$ be such that $x\lambda\in{}^\circ\PS$.
Moreover, we assume that  $M_{\T}(x\lambda)\subseteq M_{\T}(\lambda)$.
Then, we have $x\in W({\T})$.
\end{prop}

First, we prove the following lemma.

\begin{lem}
Let $\lambda\in {}^\circ\PS$ be regular and let $w\in W_\lambda$ be such
 that $w\lambda$ is dominant.
Then, we have $w\T\subsetneq\Pi_\lambda$.
\end{lem}

\proof
Assume that there is some  $\alpha\in\T$ such that  $w\alpha\not\in
\Pi_\lambda$.
Then $w{\alpha}^\vee\not\in\Pi_\lambda^\vee$.
Here, we remark that $\Pi_\lambda^\vee$ is a basis of the positive
coroot system $(\Delta^+)^\vee$.
So, there exists some $\beta,\gamma\in\Delta^+$ such that
$w{\alpha}^\vee={\beta}^\vee+{\gamma}^\vee$.
Since $w\lambda$ is dominant and regular, we have
$\langle w\lambda,{\beta}^\vee\rangle\geqslant 1$ and
$\langle w\lambda,{\gamma}^\vee\rangle\geqslant 1$.
$2\leqslant \langle
w\lambda,{\beta}^\vee+{\gamma}^\vee\rangle =\langle w\lambda,w{\alpha}^\vee\rangle=\langle \lambda,{\alpha}^\vee\rangle$.
On the other hand, $\lambda\in{}^\circ\PS$ implies $\langle
\lambda,{\alpha}^\vee\rangle=1$. This is a contradiction.\,\,\,\, $\Box$

{\it Proof of Proposition 2.3.2}

From Lemma 2.2.1, we may reduce the proposition to the case that
$\lambda$ is integral.
Put $\T_1=w\T$ and $\T_2=wx^{-1}\T$.
From Lemma 2.3.1, 
we have $\T_1\subseteq \Pi$ and
$\T_2\subseteq \Pi$.
Since $w_0 w_{\T_i}\T_i=-w_0\T_i$ holds for $i=1,2$, we have
$w_0 w_{\T_i}\T_i\subseteq \Pi$.
We put $I_1=\Ann_{U(\gggg)}(M_{\T}(\lambda))$ and $I_2=\Ann_{U(\gggg)}(M_{\T}(x\lambda))$.

From \cite{[BJ]} 4.10 Corollar, we have
$I_1=\Ann_{U(\gggg)}(M_{-w_0\T_1}(w_0
w_{\T_1} w\lambda))$
and $I_2=\Ann_{U(\gggg)}(M_{-w_0\T_2}(w_0
w_{\T_2} w\lambda))$.
Since $\langle w_0 w_{\T_i} w\lambda,{\alpha}^\vee\rangle<0 $ for
all $\alpha\in \Delta^+-w_0 w_{\T_i}\langle\T_i\rangle$, $M_{w_0
w_{\T_i}\T_i}(w_0 w_{\T_i} w\lambda)$ is irreducible.
Hence, $I_1$ and $I_2$ are primitive ideals of the same Gelfand-Kirillov
dimension.
The $\tau$-invariant of $I_1$ (respectively $I_2$) is $-w_0
\T_1$ (respectively $-w_0
\T_2$).
On the other hand, $M_{\T}(x\lambda)\subseteq M_{\T}(\lambda)$ implies
$I_1\subseteq I_2$.
Hence, we have $I_1=I_2$.
Comparing the $\tau$-invariants, we have  $-w_0
\T_1=-w_0
\T_2$.
Hence, $w\T=\T_1=\T_2=wx^{-1}\T$.
This implies $x\in W(\T)$. \,\,\,\, Q.E.D.

\subsection{Translation principle for scalar generalized Verma modules}

Next, we consider the images of scalar generalized Verma modules under certain translation functors.

For each $\gamma\in\Sigma_\T$, we put$\Delta^\gamma=\{\alpha\in\Delta^+\mid\alpha|_{\aas}=\gamma\}$. 
We prove:
\begin{lem} 
Let $\gamma\in\Sigma_\T$ and let $\beta\in \Delta^\gamma$.
If $\langle \rho_\T,\beta^\vee\rangle<0$ (resp.\ $\langle \rho_\T,\beta^\vee\rangle>0$), then there exists some $\beta^\prime\in \Delta^\gamma$ such that $\langle \rho_\T,{\beta^\prime}^\vee\rangle=\langle \rho_\T,\beta^\vee\rangle+1$ (resp.\ $\langle \rho_\T,{\beta^\prime}^\vee\rangle=\langle \rho_\T,\beta^\vee\rangle-1$).
\end{lem}
\proof  We assume that $\beta\in \Delta^\gamma$.
If $\langle \rho_\T,\beta^\vee\rangle<0$.  So, there exists some $\delta\in\T$ such that $\langle\delta,\beta^\vee\rangle<0$.  this implies that $\beta^\vee+\delta^\vee
\in\Delta^\vee$.  Hence there is some $\beta^\prime\in\Delta$ such that $\beta^\vee+\delta^\vee={\beta^\prime}^\vee$.  This $\beta^\prime$ satisfies the desirable conditions.  The remaining statement is proved in a similar way.
\,\,\,$\Box$

\begin{lem}
Let $\T\subsetneq\Pi$ and let $\mu\in\aad$ be such that $\rho_\T+\mu$ is regular integral. Let $\gamma\in\Sigma_\T$ be such that $\langle\mu,\gamma\rangle>0$.
Then for each $\beta\in\Delta^\gamma$ we have $\langle\rho_\T+\mu,\beta\rangle>0$.
\end{lem}
\proof
Put $M_\gamma=\{ \beta\in\Delta^\gamma\mid \langle\rho_\T+\mu,\beta^\vee\rangle<0\}$.  Since $\rho_\T+\mu$ is regular, we have only to show $M_\gamma=\emptyset$.
Assuming that $M_\gamma\neq\emptyset$,  we deduce a contradiction. 
We choose$\beta_0\in M_\gamma$ such that $\langle\rho_\T+\mu,\beta_0^\vee\rangle$ is maximal among the elements of $M_\gamma$.  Since $\langle\mu,\beta_0\rangle=\langle\mu,\gamma\rangle>0$, we may apply Lemma 2.4.1.  So, there exists some $\beta^\prime\in\Delta^\gamma$ such that $\langle \rho_\T+\mu,{\beta^\prime}^\vee\rangle=\langle \rho_\T+\mu,\beta_0^\vee\rangle+1$.  Since $\rho_\T+\mu$ is integral, $\langle\rho_\T+\mu,\beta_0^\vee\rangle\leqslant -1$.  If $\langle\rho_\T+\mu,\beta_0^\vee\rangle< -1$, then $\beta^\prime\in M_\gamma$.  It contradicts the choice of $\beta_0$.
If $\langle\rho_\T+\mu,\beta_0^\vee\rangle= -1$, we have $\langle \rho_\T+\mu,{\beta^\prime}^\vee\rangle=0$.  it contradicts our assumption $\rho_\T+\mu$ is regular. \,\,\,\, $\Box$

\begin{lem} 
Let $\T\subsetneq\Pi$ and let $\mu\in\aad$ be such that $\rho_\T+\mu$ is regular integral.
Then $\rho_\T+\mu$ is dominant if and only if $\langle\mu,\gamma\rangle>0$ for all $\gamma\in\Sigma^+_\T$.
\end{lem}
\proof

First, we assume that $\rho_\T+\mu$ is dominant.  We fix an arbitrary  $\gamma\in\Sigma^+_\T$.  From Lemma 2.4.1, there exists some $\beta\in\Delta^\gamma$ such that $\langle\rho_\T,\beta\rangle\geqslant 0$.
Then, we have $\langle \mu,\gamma\rangle=\langle\mu,\beta\rangle\geqslant\langle\rho_\T+\mu,\beta\rangle>0$.

Next, we assume that $\langle\mu,\gamma\rangle>0$ for all $\gamma\in\Sigma^+_\T$.  From 2.4.2, we see that $\langle\rho_\T+\mu,\beta\rangle>0$ for all $\beta\in\Delta^+-\itg\T$. On the other hand, $\langle\rho_\T+\mu,\alpha^\vee\rangle=1$ for all $\alpha\in\T$.
So, $\rho_\T+\mu$ is dominant.
   ,\,\,\,$\Box$

We put $\Sigma_\T^+(\nu)=\{\gamma\in\Sigma_\T\mid\langle\nu,\gamma\rangle>0\}$  for $\nu\in\aad$ such that $\langle\nu,\gamma\rangle\neq 0$ for all $\gamma\in\Sigma_\T$.

The following result is immediately deduced from Lemma 2.4.2.
\begin{lem}
Let $\T\subsetneq\Pi$ and let $\mu,\nu\in\aad$ and $w\in W$ be such that $\rho_\T+\mu$ is regular integral and $w(\rho_\T+\mu)=\rho_\T+\nu$. 
Then, we have $w=e$ if and only if $\Sigma_\T^+(\mu)= \Sigma_\T^+(\nu)$.
\end{lem}

We also put
\[K(\T)=\{w\in W\mid w\T\subseteq\Pi\}.\]
For $\lambda\in \hhd$ which is regular integral, we put $\Delta^+(\lambda)=\{\alpha\in\Delta\mid\langle\alpha,\lambda\rangle>0\}$

We consider the following condition (T) on $\mu,\nu\in\aad$. 
\medskip

{\bf Condition (T)}
$\rho_\T+\mu$ and $\rho_\T+\nu$ is integral and There exists some $\lambda\in\aad$ which satisfies the following (T1)-(T3).
\begin{itemize}
\item[(T1)] $\rho_\T+\lambda$ is regular integral.
\item[(T2)] We have $\langle\mu,\gamma\rangle\geqslant 0$ and
$\langle\nu,\gamma\rangle\geqslant 0$ for each $\gamma\in \Sigma_\T^+(\lambda)$.
\item[(T3)] $\mu-\nu$ is dominant with respect to $\Delta^+(\rho_\T+\lambda)$.
\end{itemize}

We fix $\mu,\nu\in\aad$ satisfying (T).
The irreducible finite-dimensional $U(\ggg)$-module $V=V_{\nu-\mu}$ with a extreme weight $\nu-\mu$ has a filtration $0=V(0)\subseteq V(1)\subseteq\cdots\subseteq V(k-1)\subseteq V(k)=V$ of $U(\pps)$-submodules such that $V(i)/V(i-1)$ \,\, ($1\leqslant i\leqslant k$) is an irreducible $U(\lls)$-modules with the highest weight $\mu_i$.  Then, there is some $1\leqslant i_0\leqslant k$ such that $\mu_{i_0}=\mu-\nu$ and $\mu_i\neq\mu-\nu$ for all $1\leqslant i\leqslant k$ such that $i\neq i_0$.
We also see that $V(i_0)/V(i_0-1)$ is a one-dimensional. 

 Let $z\in W$ be such that $z(\rho_\T+\lambda)$ is dominant.    We put $\T^\prime=z\T$,  $\nu^\prime=z\nu$and $\mu^\prime=z\mu$.
Then, from Lemma 2.3.3, we have $\T^\prime\subseteq\Pi$.  Hence $\mu^\prime\in\aaa_{\T^\prime}^\ast$ and $z(\rho_\T+\mu)=\rho_{\T^\prime}+\mu^\prime$.
We also put $\mu_i^\prime=z\mu_i$.  Twisting by $z$, $V|_{\lll_{\T^\prime}}$ decompose into the direct product of irreducible $\llll_{\T^\prime}$-modules with highest weights $\mu_1^\prime,...,\mu_k^\prime$.

The following result is more or less easy consequence of the argument of the proof of \cite{[VU]} Proposition 8.5.
\begin{lem}\,\,\, {\rm (\cite{[VI]} Lemma 4.8, \cite{[Bi]} p21 Claim, also see \cite{[M2]} Lemma 1.2.2)}
 We retain the above settings.
Assume $y\in W$.
Then, $y(\rho_{\T^\prime}+\nu^\prime)= \mu_i^\prime+\rho_{\T^\prime}+\mu^\prime$ if and only if $y(\rho_\T+\nu^\prime)=\rho_\T+\nu^\prime$ and $i=i_0$.
\end{lem}

When we regard $V=V_{\nu-\mu}$ as a $U(\pps)$-module, we write it by $V|_{\pps}$.
Since $M_\T(\rho_\T+\mu)\otimes V\cong U(\ggg)\otimes_{U(\pps)}(E_\T(\mu-\rho^\T)\otimes V|_{\pps})$ holds, we easily see the following from Lemma 2.4.5.
\begin{prop}
Let $\mu,\nu\in\aad$ be such that they satisfy the condition (T) above.

Then, we have
\[ T^{\rho_\T+\nu}_{\rho_\T+\mu}(M_\T(\rho_\T+\mu))=M_\T(\rho_\T+\nu).\]
\end{prop}

Finally, we have the following result.

\begin{prop}
Let $\lambda\in\aad$ be such that $\rho_\T+\lambda$ is regular integral and let $w\in W(\T)$.
We assume that $M_\T(\rho_\T+w\lambda)\subseteq M_\T(\rho_\T+\lambda)$.
Let $\nu\in\aad$ be such that $\rho_\T+\nu$ is integral and $\langle\nu,\gamma\rangle\geqslant 0$ for each $\gamma\in \Sigma_\T^+(\lambda)$.
Then,we have $M_\T(\rho_\T+w\nu)\subseteq M_\T(\rho_\T+\nu)$.
\end{prop}
\proof
Since $2\rho_\T$ is integral, so is $2\lambda$.  Hence for all $k\in\nat$, $\rho_\T+(2k+1)\lambda$ is regular integral.   It is easy to see $\mu=(2k+1)\lambda$ and $\nu$ satisfy the condition (T).
From the translation principle, we have $M_\T(\rho_\T+w\mu)\subseteq M_\T(\rho_\T+\mu)$.  Hence, applying proposition 2.4.7 and the exactness of the translation functor, we have the desired conclusion. \,\,\,\,$\Box$

\setcounter{section}{3}
\setcounter{subsection}{0}

\section*{\S\,\, 3.\,\,\,\,  Some results on Bruhat orderings}

\subsection{Quasi subsystems}

Let $(W_i, S_i)$ \,\,\, $(i=1,2)$ be finite coxeter systems.  We denote by $\ell_i(w)$
 the length of a reduced expression of $w\in W_i$ with respect to $S_i$.
We also denote by $\leqslant_i$ the Bruhat ordering for $(W_i, S_i)$.

\begin{dfn}
We say that $(W_2,S_2)$ is a quasi subsystem of $(W_1,S_1)$, if the following (Q1) and (Q2) hold.
\begin{itemize}
\item[(Q1)]  $W_2$ is a subgroup of $W_1$.
\item[(Q2)] For any reduced expression $w=s_1\cdots s_k$ of $w\in W_2$ in $(W_2,S_2)$, we have $\ell_1(w)=\ell_1(s_1)+\cdots+\ell_1(s_k)$.
\end{itemize}
\end{dfn}

The following lemma is easy.

\begin{lem}
Assume that $(W_2,S_2)$ is a quasi subsystem of $(W_1,S_1)$.
Then, $x\leqslant_2 y$ implies  $x\leqslant_1 y$ for all $x,y\in W_2$.
\end{lem}

We have the following lemma.

\begin{lem}
Assume that $(W_2,S_2)$ is a quasi subsystem of $(W_1,S_1)$. 
Moreover, we assume the following condition (C).
\begin{itemize}
\item[(C)]For any $x,y\in W_2$ and $s\in S_2$ such that $x\leqslant_1y$, $\ell_1(sy)<\ell_1(y)$, and $\ell_1(x)<\ell_1(sx)$, we have $sx\leqslant_1 y$.
\end{itemize}
Then,  $x\leqslant_1 y$ implies  $x\leqslant_2 y$ for all $x,y\in W_2$.
\end{lem}

\proof
Let $x,y\in W_2$ be such that $x\leqslant_1 y$.  
We show $x\leqslant_2 y$ by a double induction with respect to $\ell_2(y)$ and $\ell_2(y)-\ell_2(x)$.

Obviously we may assume $\ell_2(y)>0$.  So, we choose some $s\in S_2$ such that $\ell_2(sy)<\ell_2(y)$.

First, we assume that $\ell_2(sx)<\ell_2(x)$. We fix  reduced expressions of $s$, $sx$, and $sy$ in $(W_1, S_1)$ as follows.
\begin{align*}
s &=s_1\cdots s_k    \,\,\,\,\,\,\,\,\,    (s_1,...,s_k\in S_1), \\
sx &=t_1\cdots t_h    \,\,\,\,\,\,\,\,\,    (t_1,...,t_h\in S_1), \\
sy &=r_1\cdots r_n    \,\,\,\,\,\,\,\,\,    (r_1,...,r_n\in S_1)
\end{align*}
From (Q2), we easily see that $s_m\cdots s_kt_1\cdots t_h$ and $s_m\cdots s_k
 r_1\cdots r_n$ are reduced expressions for all $1\leqslant m\leqslant k$.
Applying \cite{[De]} Theorem 1.1, we have $s_m\cdots s_kt_1\cdots t_h\leqslant_1 s_m\cdots s_k r_1\cdots r_n$ by the induction on $m$.  So, we have $sx\leqslant_1 sy$.  Since $\ell_2(sy)<\ell_2(y)$, the induction hypothesis implies that $sx\leqslant_2 sy$. Again, applying  \cite{[De]} Theorem 1.1, we have $x\leqslant_2 y$.

Next, we assume that $\ell_2(sx)>\ell_2(x)$. From (Q2), we have $\ell_1(sx)>\ell_1(x)$.  So, we have $sx\leqslant_1y$ from (C).
Since $\ell_2(y)-\ell_2(sx)<\ell_2(y)-\ell_2(x)$, we have $sx\leqslant_2 y$ from the induction hypothesis.  Since $x\leqslant_2 sx$, we have $x\leqslant_2 y$. \,\,\,\, \hbox{$\blacksquare$}

\subsection{$\T$-useful roots}

In this subsection, we use the notation in \S1. 

Following Knapp \cite{[K]}, Howlet \cite{[H]}, and Lusztig \cite{[L]}, we consider useful roots for our purpose.

Hereafter, we fix a subset $\T$ of $\Pi$.
For $\alpha\in\Delta$, we put
\[ \Delta(\alpha)=\{\beta\in\Delta\mid\exists c\in\rel \,\,\,\, \beta|_{\aas}=c\alpha|_{\aas}\},\]
\[\Delta^+(\alpha)=\Delta(\alpha)\cap\Delta^+,\]
\[U_\alpha=\cpx \T+\cpx\alpha\subseteq\hhd.\]
Then $(U_\alpha,\Delta(\alpha),\langle\,\,\,,\,\,\,\rangle)$ is a subroot system of $(\hhh^\ast, \Delta,\langle\,\,\,,\,\,\,\rangle)$.
The set of simple roots for $\Delta^+(\alpha)$ is denoted by $\Pi(\alpha)$.
$\alpha|_{\aas}=0$ if and only if $\T=\Pi(\alpha)$. 
For $\alpha\in\Delta^+$, we denote by $W_{\T}(\alpha)$ the Weyl group of $(\hhh^\ast, \Delta(\alpha))$.  
Clearly, $W_{\T}\subseteq W_{\T}(\alpha)\subseteq W$.
We denote by $w^\alpha$ the longest element of $W_{\T}(\alpha)$.
We put as follows.
$$\sigma_\alpha=w^\alpha w_{\T}.$$
$\alpha|_{\aas}=0$ if and only if $\sigma_\alpha=1$.
If $\alpha\in\Delta$ is orthogonal to all the elements in $\T$, then we can easily see $\alpha$ is $\T$-reduced and $s_\alpha=\sigma_\alpha$.
\begin{dfn}
\begin{itemize}
\item[{\rm(1)}]{\rm We call $\alpha\in\Delta$ $\T$-{\it useful} if the order of ${\sigma_\alpha}$ is two. We denote by ${}^u\Delta_\T$ the set of the useful $\T$-roots. We also put ${}^u\Delta_\T^+={}^u\Delta_\T\cap\Delta^+$.}
\item[{\rm (2)}]{\rm If $\alpha|_{\aas}\neq 0$, then $\Pi(\alpha)$ is written as $S\cup\{\tilde{\alpha}\}$.
If $\alpha\in\Delta$ satisfies $\alpha|_{\aas}\neq 0$ and $\alpha=\tilde{\alpha}$, then we call $\alpha$ $\T$-{\it reduced.}  We put}
$$ {}^{ru}\Delta_\T^+=\{\alpha\in {}^u\Delta_\T^+\mid \hbox{$\alpha$ is $\T$-reduced.}\}$$
\end{itemize}
\end{dfn}

We easily see:
\begin{lem}
Let $\alpha\in\Delta^+$ be $\T$-reduced.  We denote by $\Delta(\alpha)_0$ be the irreducible component of $\Delta(\alpha)$ containing $\alpha$. We put $\Pi(\alpha)_0=\Pi(\alpha)\cap\Delta(\alpha)_0$.
\begin{itemize}
\item[(1)] If $\Delta(\alpha)_0$ is not of the type ADE, then we have $\alpha\in{}^{ur}\Delta_\T^+$.
\item[(2)] If $\Delta(\alpha)_0$ is of the type ${\rm D}_{2n}$ \,\,$(n\geqslant 2)$, ${\rm E}_7$, or ${\rm E}_8$, then we have $\alpha\in{}^{ur}\Delta_\T^+$.
\item[(3)] If $\Delta(\alpha)_0$ is of the type ${\rm A}_{2n}$ \,\,$(n\geqslant 1)$, then we have $\alpha\not\in{}^{ur}\Delta_\T^+$.
\item[(4)] We assume that $\Delta(\alpha)_0$ is of the type ${\rm A}_{2n+1}$ \,\,$(n\geqslant 0)$.  We number the elements of $\Pi(\alpha)_0$ as follows.
$$ \Pi(\alpha)_0=\{\beta_1,....,\beta_{2n+1}\}.$$
We choose  the above numbering so that $\langle\beta_i,\beta_{i+1}^{\vee}\rangle=-1$ for $1\leqslant i\leqslant 2n$.  Then $\alpha\in{}^{ur}\Delta_\T^+$ if and only if $\alpha=\beta_n$.
\item[(5)] We assume that $\Delta(\alpha)_0$ is of the type ${\rm D}_{2n+1}$ \,\,$(n\geqslant 2)$.
 We number the elements of $\Pi(\alpha)_0$ as follows.
$$ \Pi(\alpha)_0=\{\beta_1,....,\beta_{2n+1}\}.$$
We choose  the above numbering so that $\langle\beta_i,\beta_{i+1}^{\vee}\rangle=-1$ for $1\leqslant i\leqslant 2n-1$ and  $\langle\beta_{2n-1},\beta_{2n+1}^{\vee}\rangle=-1$ .  Then $\alpha\in{}^{ur}\Delta_\T^+$ if and only if $\alpha\not\in\{\beta_{2n},\beta_{2n+1}\}$.
\item[(6)] We assume that $\Delta(\alpha)_0$ is of the type ${\rm E}_{6}$.
 We number the elements of $\Pi(\alpha)_0$ as follows.
$$ \Pi(\alpha)_0=\{\beta_1,....,\beta_{6}\}.$$
We choose  the above numbering so that $\langle\beta_i,\beta_{i+1}^{\vee}\rangle=-1$ for $1\leqslant i\leqslant 4$ and  $\langle\beta_{3},\beta_{6}^{\vee}\rangle=-1$ .  Then $\alpha\in{}^{ur}\Delta_\T^+$ if and only if $\alpha\in\{\beta_{3},\beta_{6}\}$.
\end{itemize}
\end{lem}
For $\alpha\in {}^{ru}\Delta_\T$, we put 
\[ V_\alpha=\{\lambda\in\ads\mid\langle\lambda,\alpha\rangle=0\},\]
$$\hat{\alpha}=\alpha|_{\aas} \in\ads.$$

We easily see:
\begin{lem}
Let $\alpha\in {}^{ru}\Delta_\T^+$.  Then, we have

\begin{itemize}
\item[(1)] $\sigma_\alpha$ preserves $\ads$.
\item[(2)] $\sigma_\alpha\in W(\T)$.  In particular, $\sigma_\alpha\rho_{\T}=\rho_{\T}$. 
\item[(3)]$\sigma_{\alpha}\hat{\alpha}=-\hat{\alpha}$.
\item[(4)] $\sigma_\alpha|_{\ads}$ is the reflection with respect to $V_\alpha$.
\end{itemize}
\end{lem}
We denote by $W(\T)^\prime$ the subgroup of $W$ generated by $\{\sigma_\alpha\mid \alpha\in{}^{ru}\Delta^+_\T\}$.
We put ${}^{u}\Sigma_\T=\{\alpha|_{\aas}\in\aad\mid \alpha\in {}^u\Delta_\T\}$. 
${}^u\Sigma_\T$ is a (not necessarily reduced) root system. 
We also put ${}^{ru}\Sigma_\T^+=\{\alpha|_{\aas}\in\aad\mid \alpha\in {}^{ru}\Delta_\T^+\}$ and ${}^{ru}\Sigma_\T={}^{ru}\Sigma_\T^+\cup-{}^{ru}\Sigma_\T^+$.
${}^{ru}\Sigma_\T$ is a reduced root system and ${}^{ru}\Sigma^+_\T$ is a positive system.  
We denote by ${}^u\Pi_\T$ the simple system for ${}^{ru}\Sigma^+_\T$.
${}^u\Pi_\T$ is also a basis of ${}^{u}\Sigma_\T$.
For $\alpha\in {}^{ru}\Delta_\T^+$, $\sigma_\alpha$ depends only on $\alpha|_{\aas}$.
So, sometimes we write $\sigma_{\alpha|_{\aas}}$ for $\sigma_\alpha$.
We put $S(\T)=\{\sigma_\gamma\mid \gamma\in {}^u\Pi_\T\}$

\begin{thm} {\rm (Howlet \cite{[H]} Theorem 6, Lusztig \cite{[L]} \S 5)}

(1)  \,\,\, $W(\T)^\prime\subseteq W(\T)$.

(2) \,\,\, For $\alpha\in{}^u\Delta^+_\T$, $\sigma_\alpha(\aas^\ast)=\aas^\ast$.
Moreover, $\sigma_\alpha|_{\aas^\ast}$ is the reflection  with respect to $\alpha|_{\aas}$ and $\sigma_\alpha\rho_\T=\rho_\T$.

(3) \,\,\, We define $\iota : W(\T)^\prime\rightarrow\hbox{GL}(\aas^\ast)$ by $\iota(x)=x|_{\aas^\ast}$.  Then $\iota$ is an injective group homomorphism.

(4) \,\,\, $\iota(W(\T)^\prime)$ is the reflection group for the root system ${}^{ru}\Sigma_\T$.  Hence $(W(\T)^\prime, S(\T))$ is a Coxeter system.
\end{thm}

We denote by $\leqslant_\T$ the Bruhat ordering for $(W(\T)^\prime, S(\T))$.

\subsection{Normal parabolic subalgebras}

\begin{dfn}
We call $\T\subsetneq\Pi$ normal, if $\Pi-\T\subseteq {}^u\Delta_\T^+$.
We call a standard parabolic subalgebra $\pps$ normal, if $\T$ is normal.
A parabolic subalgebra is called normal, if it is conjugate to a normal standard parabolic subalgebra by an inner automorphism.
\end{dfn}

We describe the list of the  normal parabolic
subalgebras of   classical Lie algebras.

(1) \,\, Let $\gggg=\gggg\llll(n,\cpx)$ (the case of $\gggg=\sss\llll(n,\cpx)$
    is similar) and let $k$ be a positive integer dividing $n$.
 We consider the following parabolic subalgebras.

$\ppp(A_{n-1,k})$  : a parabolic subalgebra of $\gggg$ whose Levi part is
isomorphic to

$$\overbrace{\gggg\llll(k,\cpx)\oplus\cdots\oplus\gggg\llll(k,\cpx)}^{{n/k}}.$$

(2) \,\, Let $\gggg$ be a complex simple Lie algebra of the type $X_n$.
Here, $X$ means one of $B$, $C$, and $D$.
Let $k$ and $\ell$ be  positive integers such that $k$ divides $n-\ell$.
If $X=D$, then we assume that $\ell\neq 1$.

 We consider the following parabolic subalgebras.

$\ppp(X_{n,k,\ell})$  :  a parabolic subalgebra of $\gggg$ whose Levi part is
isomorphic  to

$$\overbrace{\gggg\llll(k,\cpx)\oplus\cdots\oplus\gggg\llll(k,\cpx)}^{{(n-\ell)}/k}\oplus
X_{\ell}.$$

Here, $X_{\ell}$ means that the complex simple Lie algebra of the
type $X_{\ell}$.  Namely $B_\ell={\mathfrak so}(2\ell+1,\cpx)$, $C_n={\mathfrak sp}(\ell,\cpx)$, and $D_n={\mathfrak so}(2\ell,\cpx)$.  $X_0$ means the zero Lie algebra.

From lemma 3.2.2, we easily see:
\begin{prop}
\begin{itemize}
\item[(1)] $\ppp(A_{n-1,k})$ is normal.
 Conversely any normal parabolic subalgebra of  is conjugate to $\ppp(A_{n,k})$
 for some $k$.
\item[(2)] 
$\ppp(X_{n,k,\ell})$ is normal   $X=D$, $\ell=0$, and
$k$ is an odd number greater than $1$.
Any normal parabolic subalgebra is conjugate to one of such
$\ppp(X_{n,k,\ell})$s by an inner automorphism.
\end{itemize}
\end{prop}

For exceptional simple Lie algebras, we have the following results.
We describe $\T\subsetneq\Pi$ by a marked Dynkin diagram.
We write "\textcircled{$\bullet$}" the vertices corresponding to elements in $\T$.
If $\T$ is the empty set, it is obviously normal.  So, we consider $\emptyset\neq\T\subsetneq\Pi$.
\begin{prop}

\begin{itemize}
\item[(1)] Assume that $\ggg$ is of type ${\rm G}_2$.
Then any subset of $\Pi$ is normal.
\item[(2)] Assume that $\ggg$ is of type ${\rm F}_4$.
If $\card\T=3$, $\T\subsetneq\Pi$ is normal.
The list of the other nonempty normal subsets of $\Pi$ is as follows.
\begin{itemize}
\item[${\rm F}_{4,12}$]
$$\bigcirc-\bigcirc\Leftarrow\hbox{\textcircled{$\bullet$}}-\hbox{\textcircled{$\bullet$}}$$
\item[${\rm F}_{4,14}$]
$$\bigcirc-\hbox{\textcircled{$\bullet$}}\Leftarrow\hbox{\textcircled{$\bullet$}}-\bigcirc$$
\item[${\rm F}_{4,34}$]
$$\hbox{\textcircled{$\bullet$}}-\hbox{\textcircled{$\bullet$}}\Leftarrow\bigcirc-\bigcirc$$
\end{itemize}
\item[(3)] Assume that $\ggg$ is of type ${\rm E}_6$.
Then the list of the nonempty normal subsets of $\Pi$ is as follows.
\begin{itemize}
\item[${\rm E}_{6,3}$]
$$\begin{array}{c}
\textcircled{$\bullet$}-\textcircled{$\bullet$}-\bigcirc-\textcircled{$\bullet$}-\textcircled{$\bullet$}\\|\\ \textcircled{$\bullet$}
\end{array}$$
\item[${\rm E}_{6,6}$]
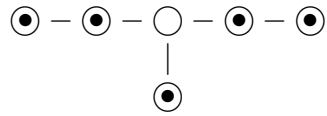
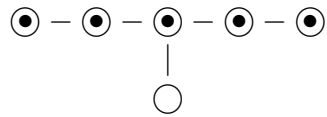
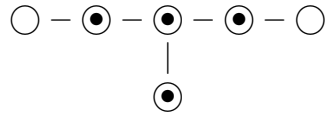
$$\begin{array}{c}
\textcircled{$\bullet$}-\textcircled{$\bullet$}-\textcircled{$\bullet$}-\textcircled{$\bullet$}-\textcircled{$\bullet$}\\|\\ \bigcirc
\end{array}$$
\item[${\rm E}_{6,15}$]
$$\begin{array}{c}
\bigcirc-\textcircled{$\bullet$}-\textcircled{$\bullet$}-\textcircled{$\bullet$}-\bigcirc\\|\\ \textcircled{$\bullet$}
\end{array}$$
\end{itemize}
\item[(4)]Assume that $\ggg$ is of type ${\rm E}_7$.
If $\card\T=6$, $\T\subsetneq\Pi$ is normal.
The list of the other nonempty normal subsets of $\Pi$ is as follows.
\begin{itemize}
\item[${\rm E}_{7,27}$]
$$\begin{array}{cl}
\bigcirc-\textcircled{$\bullet$}-\textcircled{$\bullet$}-\textcircled{$\bullet$}-\bigcirc&-\textcircled{$\bullet$}\\|&\\ \textcircled{$\bullet$}&
\end{array}$$
\item[${\rm E}_{7,67}$]
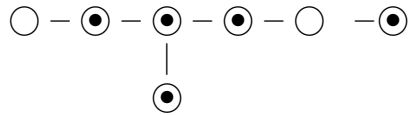
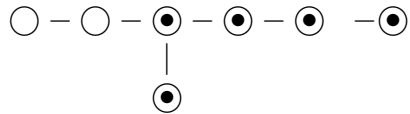
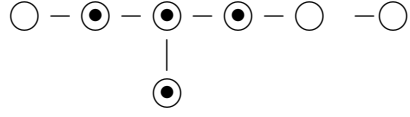
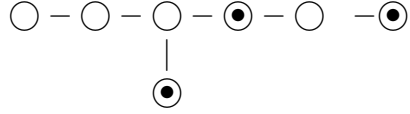
$$\begin{array}{cl}
\bigcirc-\bigcirc-\textcircled{$\bullet$}-\textcircled{$\bullet$}-\textcircled{$\bullet$}&-\textcircled{$\bullet$}\\|&\\ \textcircled{$\bullet$}&
\end{array}$$
\item[${\rm E}_{7,127}$]
$$\begin{array}{cl}
\bigcirc-\textcircled{$\bullet$}-\textcircled{$\bullet$}-\textcircled{$\bullet$}-\bigcirc&-\bigcirc\\|&\\ \textcircled{$\bullet$}&
\end{array}$$
\item[${\rm E}_{7,2467}$]
$$\begin{array}{cl}
\bigcirc-\bigcirc-\bigcirc-\textcircled{$\bullet$}-\bigcirc&-\textcircled{$\bullet$}\\|&\\ \textcircled{$\bullet$}&
\end{array}$$
\end{itemize}
\item[(5)]Assume that $\ggg$ is of type ${\rm E}_8$.
If $\card\T=7$, $\T\subsetneq\Pi$ is normal.
The list of the other nonempty normal subsets of $\Pi$ is as follows.
\begin{itemize}
\item[${\rm E}_{8,12}$]
$$\begin{array}{cl}
\textcircled{$\bullet$}-\textcircled{$\bullet$}-\textcircled{$\bullet$}-\textcircled{$\bullet$}-\textcircled{$\bullet$}&-\bigcirc-\bigcirc\\|&\\ \textcircled{$\bullet$}&
\end{array}$$
\item[${\rm E}_{8,18}$]
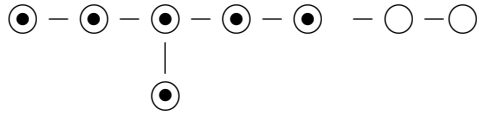
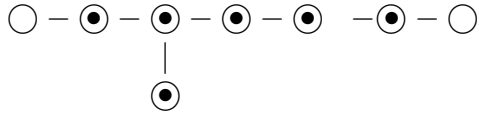
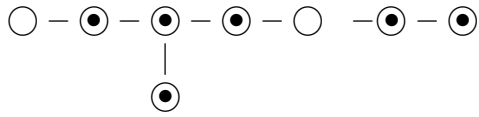
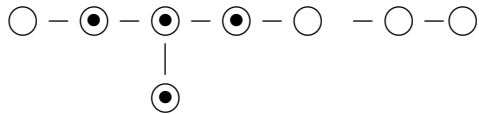
$$\begin{array}{cl}
\bigcirc-\textcircled{$\bullet$}-\textcircled{$\bullet$}-\textcircled{$\bullet$}-\textcircled{$\bullet$}&-\textcircled{$\bullet$}-\bigcirc\\|&\\ \textcircled{$\bullet$}&
\end{array}$$
\item[${\rm E}_{8,38}$]
$$\begin{array}{cl}
\bigcirc-\textcircled{$\bullet$}-\textcircled{$\bullet$}-\textcircled{$\bullet$}-\bigcirc&-\textcircled{$\bullet$}-\textcircled{$\bullet$}\\|&\\ \textcircled{$\bullet$}&
\end{array}$$
\item[${\rm E}_{8,1238}$]
$$\begin{array}{cl}
\bigcirc-\textcircled{$\bullet$}-\textcircled{$\bullet$}-\textcircled{$\bullet$}-\bigcirc&-\bigcirc-\bigcirc\\|&\\ \textcircled{$\bullet$}&
\end{array}$$
\end{itemize}
\end{itemize}
\end{prop}

We give some characterization of normality.

\begin{prop}
For $\T\subsetneq\Pi$. the following conditions are equivalent.
\begin{itemize}
\item[(1)]  $\T\subsetneq\Pi$ is normal.
\item[(2)]  $K(\T)=W(\T)^\prime$.
\item[(3)]  $K(\T)=W(\T)$.
\item[(4)]  ${}^u\Sigma_\T=\Sigma_\T$
\end{itemize}
\end{prop}

\proof
First, we assume (1) .
Then, using Proposition 3.3.2 and 3.3.3, we obtain (2) and (4) via case-by-case analysis.
(2) obviously implies (3).
Next, we assume (3).
For $\alpha\in \Pi-\T$, we easily see $\sigma_\alpha^2(\Pi)\subseteq\Delta^+$.
Hence $\sigma_\alpha$ is an involution.  This means that $\alpha\in{}^u\Delta_\T^+$.  So, we have (1).  
(4) is clearly stronger than (1).
\,\,\,\,$\Box$
\begin{cor}
If $\T\subsetneq\Pi$ is normal, then $W(\T)^\prime=W(\T)$.
\end{cor}

Since $\Delta^+\cap(-w\Delta^+)=\{\alpha\in\Delta^+\mid\alpha|_{\aas}\in\Sigma_\T^+\cap (-w\Sigma_\T^+)\}$ for each $w\in W(\T)$,
we easily see the following lemma.
\begin{lem} We assume that $\T\subsetneq\Pi$ is normal.
Then for each $w\in W(\T)$, we have 
$$\Delta^+\cap(-w\Delta^+)=\bigcup
	_{\gamma\in{}^{ru}\Sigma_\T^+\cap(-w{}^{ru}\Sigma_\T^+)}\{\alpha\in\Delta^+\mid \exists c>0 \,\,\, \alpha|_{\aas}=c\gamma\}.$$
\end{lem}
Hence, we have the following result.
\begin{prop}
If $\T\subsetneq\Pi$ is normal, then $(W(\T)^\prime, S(\T))$ is a quasi subsystem of $(W, S)$.
\end{prop}
As a corollary of Proposition 3.3.4, we easily have:
\begin{cor}
$\T\subsetneq\Pi$ is normal if and only if any two parabolic subalgebras with the Levi part $\llll_\T$ are conjugate under an inner automorphism of $\ggg$.
\end{cor}

\subsection{Comparison of Bruhat orderings}

In this subsection, we use the notation in \S1. 
\begin{dfn}
{\rm
We call $\T\subsetneq\Pi$ {\it seminormal}, if there exists some $\Psi$ such that $\T\subseteq\Psi\subseteq\Pi$ and ${}^u\Pi_\T=\{\alpha|_{\aas}\mid \alpha\in\Psi-\T\}$. 
}
\end{dfn}
 So, $S(\T)=\{\sigma_\alpha\mid \alpha\in\T-\Psi\}$.

$\T\subsetneq\Pi$ is seminormal if and only if there is a $\alpha\in\Pi\cap{}^{ru}\Delta^+$ such that $\alpha|_{\aas}=\gamma$ for each $\gamma\in {}^u\Pi_\T$.

We immediately see the following result from Proposition 3.3.7.
\begin{cor}
If $\T\subsetneq\Pi$ is seminormal, then $(W(\T)^\prime, S(\T))$ is a quasi subsystem of $(W, S)$.
\end{cor}
We fix a  connected  complex reductive Lie group $G$ whose Lie algebra is $\ggg$.  For $\T\subsetneq\Pi$, we denote by $P_\T$ (resp.\ $H$) the parabolic subgroup (resp.\ the Cartan subgroup) of $G$ corresponding to $\pps$ (resp.\ $\hhh$).  We denote by $N_G(H)$ the normalizer of $H$ in $G$.   Since the Weyl group $W$ is identified with the quotient group $N_G(H)/H$, for each $w\in W$ we can fix a representative in $N_G(H)$.  We denote the representative by the same letter "$w$".
 
For $x\in W$, we put $U_x=P_\T x/P_\T$.  Namely, $U_x$ is a $P_\T$-orbit in $G/P_\T$ through $x/P_\T\in G/P_\T$.
We denote by $\overline{U}_x$ the closure of $U_x$ in $G/P_\T$.
If $w\in W(\T)$, then $\ell(ws_\alpha)>\ell(w)$ for all $\alpha\in \T$.
Hence, we have
\begin{lem}
\begin{itemize}
\item[(1)] For $w\in W(\T)$, we have $\dim U_w=\ell(w)$.
\item[(2)] For $x, y\in W(\T)$, $x\leqslant y$ if and only if $\overline{U}_x\subseteq\overline{U}_y$.
\end{itemize}
\end{lem}
Next we show,
\begin{lem}
Assume that $\T\subsetneq\Pi$ is seminormal. We choose  $\T\subseteq\Psi\subseteq\Pi$ as in 3.4.1.
Fix $x\in W(\T)^\prime$.  Let $\alpha\in \Psi-\T$ be such that $\ell(\sigma_\alpha x)<\ell(x)$.
Then we have $\overline{U}_x=P_{\T\cup\{\alpha\}}\overline{U}_x=P_{\T\cup\{\alpha\}}\overline{U}_{\sigma_\alpha x}$.
\end{lem}
\proof
 We may choose a reduced expression $x=\sigma_{\alpha_1}\cdots\sigma_{\alpha_k}$ such that $\alpha_1=\alpha$.
We consider a contraction map as follows.
$$ F : P_{\T\cup\{\alpha_1\}}\times_{P_\T}P_{\T\cup\{\alpha_2\}}\times_{P_\T}\cdots\times_{P_\T}P_{\T\cup\{\alpha_k\}}/P_\T \rightarrow G/P_\T.$$
We easily see :
\begin{itemize}
\item[(a)] ${\rm Image}(F)$ is an irreducible  Zariski closed set in $G/P_\T$.
\item[(b)] $\dim \overline{U}_x=\ell(x)=\dim P_{\T\cup\{\alpha_1\}}\times_{P_\T}\cdots\times_{P_\T}P_{\T\cup\{\alpha_k\}}/P_\T  $.
\item[(c)] $\overline{U}_x\subseteq {\rm Image}(F)$.
\end{itemize}
Hence, we have $\overline{U}_x= {\rm Image}(F)$.
So, we have the lemma immediately.
\,\,\,\, $\Box$

The following result is the main result of this section.

\begin{thm}
Let $\T\subsetneq\Pi$ be seminormal.
For $x,y\in W(\T)^\prime$, $x\leqslant y$ if and only if $x\leqslant_\T y$.
\end{thm}

\proof
We choose  $\T\subseteq\Psi\subseteq\Pi$ as in 3.4.1.
From Lemma 3,1.2, Lemma 3.1.3, and Corollary 3.4.2, we have only to show the condition (C) in the statement of Lemma 3.1.3 holds for $(W(\T)^\prime, S(\T))$. 
So we choose $x,y\in W(\T)^\prime$ and $\alpha\in \Psi-\T$ such that $x\leqslant y$, $\ell(\sigma_\alpha y)<\ell(y)$, and $\ell(\sigma_\alpha x)>\ell(x)$.
From $x\leqslant y$, we have $\overline{U}_x\subseteq\overline{U}_y$ by Lemma 3.4.3 (2). Hence $P_{\T\cup\{\alpha\}}\overline{U}_x\subseteq P_{\T\cup\{\alpha\}}\overline{U}_y$.
From Lemma3.4.4, we have $\overline{U}_y=P_{\T\cup\{\alpha\}}\overline{U}_y$ and $\overline{U}_{\sigma_\alpha x}=P_{\T\cup\{\alpha\}}\overline{U}_x$.
So, we have $\overline{U}_{\sigma_\alpha x}\subseteq\overline{U}_y$.
this means that $ \sigma_\alpha x\leqslant y$.
Hence, the condition (C) holds for $\T$. \,\,\,\,\, Q.E.D.

\setcounter{section}{4}
\setcounter{subsection}{0}

\section*{\S\,\, 4.\,\,\,\, Elementary homomorphisms }
\subsection{Elementary homomorphisms}

We fix a subset $\T$ of $\Pi$ and
$\alpha\in{}^{ru}\Delta_\T^+$.
We define
\begin{gather*}
\gggg(\alpha) =\hhh+\sum_{\beta\in\Delta(\alpha)}\gggg_\beta, \,\,\,\,\,\,\,\,\,
 \ppp_\T(\alpha) =\gggg(\alpha)\cap\ppp_{\T}.
\end{gather*}
Then, $\gggg(\alpha)$ is a reductive Lie subalgebra of $\gggg$ whose root system is $\Delta(\alpha)$ and  $\pps(\alpha)$ is a maximal parabolic subalgebra of  $\gggg(\alpha)$.

We denote by $\omega_\alpha\in\ads\subseteq\hhd$ the fundamental weight
for $\alpha$ with respect to the basis $\Pi(\alpha)=\T\cup\{\alpha\}$.
Namely $\omega_\alpha$ satisfies that $\langle\omega_\alpha,
\beta\rangle=0$ for $\beta\in\T$, $\langle \beta,{\alpha}^\vee\rangle=1$,
and $\omega_\alpha|_{\hhh\cap\ccc(\gggg(\alpha))}=0$.
Here, $\ccc(\gggg(\alpha))$ is the center of $\gggg(\alpha)$. 
We see that there is some positive real number $a$ such that
$\omega_\alpha=a\alpha|_{\aas}$, since
$\alpha|_{\hhh\cap\ccc(\gggg(\alpha))}=0 $.
Hence, we have $V_\alpha=\{\lambda\in\ads\mid \langle \lambda,\omega_\alpha\rangle=0\}$.

Put
$\rho(\alpha)=\frac{1}{2}\sum_{\beta\in\Delta^+(\alpha)}\beta$,
For $\nu\in\ads$, we denote by $\cpx_\nu$ the one-dimensional $U(\pps(\alpha))$-module corresponding to $\nu$.
For $\nu\in\ads$ we define a generalized Verma module for $\gggg(\alpha)$ as follows.
\begin{align*}
M_{\T}^{\gggg(\alpha)}(\rho_{\T}+\nu)=U(\gggg(\alpha))\otimes_{U(\pps(\alpha))}\cpx_{\nu-\rho(\alpha)}.
\end{align*}

Then, we have:
\begin{thm} (\cite{[M]})
Let $\nu$ be an arbitrary element in $V_\alpha$, let $c$ be either
 $1$ or $\frac{1}{2}$.
Assume that $M_{\T}^{\gggg(\alpha)}(\rho_{\T}-c\omega_\alpha)\subseteq M_{\T}^{\gggg(\alpha)}(\rho_{\T}+c\omega_\alpha)$.
Then, we have
$M_{\T}(\rho_{\T}+\nu-(c+n)\omega_\alpha)\subseteq M_{\T}(\rho_{\T}+\nu+(c+n)\omega_\alpha)$ for all $n\in\nat$.
\end{thm}

 We call the above homomorphism of
 $M_{\T}(\rho_{\T}+\nu-(c+n)\omega_\alpha)$ into
 $M_{\T}(\rho_{\T}+\nu+(c+n)\omega_\alpha)$ {\it an elementary homomorphism}. 
In \cite{[M]}, homomorphisms between scalar generalized Verma modules associated with a maximal parabolic subalgebra are classified.  So, elementary homomorphisms are understood.

The following conjecture is propsed in \cite{[M]} as a working hypothesis.

\begin{conj}
 An arbitrary nontrivial homomorphism between  scalar
generalized Verma modules is a composition of elementary homomorphisms. 
\end{conj}

The conjecture in the case of the Verma modules is nothing but the result
of Bernstein-Gelfand-Gelfand (\cite{[BGG]}).
I do not know a counterexample for the above working hypothesis and we obtain partial affirmative results in this article. 
A weaker version is :
\begin{conj}
Let $\T\subseteq\Pi$ be normal and let $\mu,\nu\in\aad$ be such that $\rho_\T+\mu$ and $\rho_\T+\nu$ are regular integral.
If $M_\T(\rho_\T+\nu)\subseteq M_\T(\rho_\T+\mu)$, then it is a composition of elementary homomorphisms.
\end{conj}
Later, we show that the conjecture is affirmative for strictly normal case (see \S5) and exceptional Lie algebras (see \S5 and \S6).

For example, I do not know whether an homomorphism of the form $M_\T(\rho_\T+\sigma_\alpha\mu)\subseteq M_\T(\rho_\T+\mu)$ \,\, $(\mu\in\aad)$ is always elementary.
We have a weak result.

\begin{prop}

We fix $\mu\in \aad$ such that $M_\T(\rho_\T+\sigma_\alpha\mu)\subseteq M_\T(\rho_\T+\mu)$ and $\rho_\T+\mu$ is regular and integral.
If $\{\beta\in\Sigma_\T- \rel\alpha|_{\aas}\mid\langle\mu,\beta\rangle>0\}=\{\beta\in\Sigma_\T-\rel\alpha|_{\aas}\mid\langle\sigma_\alpha\mu,\beta\rangle>0\}$, then $M_\T(\rho_\T+\sigma_\alpha\mu)\subseteq M_\T(\rho_\T+\mu)$ is an elementary homomorphism.
\end{prop}

\proof

Put $\nu_0=\mu-\langle\mu,\alpha^\vee\rangle\omega_\alpha$.   Then $\nu_0\in V_\alpha$.  Since $M_\T(\rho_\T+\sigma_\alpha\mu)\subseteq M_\T(\rho_\T+\mu)$, we have $\mu-\sigma_\alpha\mu=2\langle\mu,\alpha^\vee\rangle\omega_\alpha\in\Q^+$.
Hence, $2\langle\mu,\alpha^\vee\rangle\omega_\alpha$ is integral.  So, we can write  $\langle\mu,\alpha^\vee\rangle=c+n_0$.  Here, $c$ is either
 $1$ or $\frac{1}{2}$ and $n_0$ is a positive integer. Put $\kappa=2(\mu+\sigma_\alpha\mu)$
Since $2\rho_\T$ and $\rho_\T+\mu$ are integral, so is $\kappa$.
Moreover, we have $\kappa\in V_\alpha$ and $\langle\kappa,\beta\rangle>0$ for all $\beta\in\Sigma_\T- \rel\alpha|_{\aas}$ such that $\langle\mu,\beta\rangle>0$.
From the translation principle, we have $M_{\T}(\rho_{\T}+(\nu_0+m\kappa)-(c+n_0)\omega_\alpha)\subseteq M_{\T}(\rho_{\T}+(\nu_0+m\kappa)+(c+n_0)\omega_\alpha)$ for all $m\in\nat$.
Hence $\{a\in{\mathbb C}\mid M_{\T}(\rho_{\T}+(\nu_0+a\kappa)-(c+n_0)\omega_\alpha)\subseteq M_{\T}(\rho_{\T}+(\nu_0+a\kappa)+(c+n_0)\omega_\alpha)\}$ is Zariski dense in ${\mathbb C}$.  So, we can prove  $M_{\T}(\rho_{\T}+(\nu_0+a\kappa)-(c+n_0)\omega_\alpha)\subseteq M_{\T}(\rho_{\T}+(\nu_0+a\kappa)+(c+n_0)\omega_\alpha)$ for all $a\in\cpx$ in the same way as \cite{[L]} Lemma 5.4.
If $a\in\cpx$ is generic, then the integral toot system for $(\rho_{\T}+(\nu_0+a\kappa)-(c+n_0)\omega_\alpha$ is $\Delta(\alpha)$.
Hence, Lemma 2.2.1 implies that $M_{\T}^{\gggg(\alpha)}(\rho_{\T}-(c+n_0)\omega_\alpha)\subseteq M_{\T}^{\gggg(\alpha)}(\rho_{\T}+(c+n_0)\omega_\alpha)$.  Applying \cite{[M]} Lemma 2.2.6, we have $M_{\T}^{\gggg(\alpha)}(\rho_{\T}-c\omega_\alpha)\subseteq M_{\T}^{\gggg(\alpha)}(\rho_{\T}+c\omega_\alpha)$.   \,\,\,\,$\Box$

\subsection{$\T$-excellent roots}

We retain the notations in 4.1.

\begin{dfn}
{\rm 
\begin{itemize}
\item[(1)] We call $\alpha\in{}^{ru}\Delta=\T^+$ $\T$-excellent if $\sigma_\alpha$ is a Duflo involution (\cite{[D]} cf.\ \cite{[JP]}) in $W(\alpha)$.
\item[(2)] We put ${}^e\Delta_\T^+=\{\alpha\in{}^{ru}\Delta_\T^+\mid \text{$\alpha$ is $\T$-excellent}\}$.
\item[(3)]We put ${}^e\Sigma_\T^+=\{\alpha|_{\aas}\in\aad\mid \alpha\in{}^e\Delta_\T^+\}$ and ${}^e\Sigma_\T={}^e\Sigma_\T^+\cup(-{}^e\Sigma_\T^+)$.
\item[(4)] We denote by ${}^eW(\T)$ the subgroup of $W(\T)^\prime$ generated by $\{\sigma_\alpha\mid \alpha\in {}^e\Delta_\T^+\}$.
\item[(5)] For $\alpha\in {}^{ru}\Delta_\T^+$, we put $c_\alpha=1$ (resp.\ $c_\alpha=\frac{1}{2}$) if $\rho_\T$ is integral (resp.\ not integral) with respect to $\Delta(\alpha)$.  Then, $\rho_\T+(c_\alpha+n)\omega_\alpha$ is integral with respect to $\Delta(\alpha)$ for all $n\in\itg$. 
\end{itemize}
}
\end{dfn}

We have
\begin{prop} Let $\alpha\in{}^e\Delta_\T^+$ and let $\mu\in\aad$ be such that $\rho_\T+\mu$ is integral and $\langle\mu,\alpha\rangle>0$.
Then, we have an elementary homomorphism $M_\T(\rho_\T+\sigma_\alpha\mu)\subseteq M_\T(\rho_\T+\mu)$.
\end{prop}
\proof
Put $\nu_0=\mu-\langle\mu,\alpha^\vee\rangle\omega_\alpha$.   Then $\nu_0\in V_\alpha$.  Since $\rho_\T+\mu$ is integral, we have $\langle\rho_\T+\mu,\alpha^\vee\rangle\in\itg$. From he definition of $c_\alpha$, we have $\langle\rho_\T,\alpha^\vee\rangle\in c_\alpha+\itg$.  Hence, we can write $\mu=\nu_0+(c_\alpha+n)\omega_\alpha$ for some $n\in\nat$.  So, from $\alpha\in{}^e\Delta_\T^+$,  Theorem 4.1.1, Proposition 2.4.7, and \cite{[M2]} Proposition 2.1.2, we have the proposition.\,\,\,\, $\Box$

For a simple Lie algebra of the type A, every involution is a Duflo involution (\cite{[Du]}).  hence, we have:
\begin{cor}
If $\ggg$ is a simple Lie algebra of the type A, we have ${}^{ru}\Delta_\T^+={}^{e}\Delta_\T^+$ for all $\T\subsetneq\Pi$.
\end{cor}

\setcounter{section}{5}
\setcounter{subsection}{0}

\section*{\S\,\, 5.\,\,\,\, Strictly normal case }

\subsection{Strictly normal subset of $\Pi$}

\begin{dfn}
{\rm
We call $\T\subsetneq\Pi$ {\it strictly normal}, if $\T$ is normal and ${}^e\Delta_\T^+={}^{ru}\Delta_\T^+$. A standard parabolic subalgebra $\pps$ is called strictly normal when $\T$ is strictly normal.
}
\end{dfn}

Before stating the main result, we prove the following lemma.

\begin{lem}
Let $\T\subsetneq\Pi$ be normal and let $\mu\in\aad$ be such that $\rho_\T+\mu$ is integral.  Then, $\mu$ is integral with respect to ${}^{ru}\Sigma_\T$.
\end{lem}
\proof
Since $\rho_\T+\mu$ is integral, we have $w(\rho_\T+\mu)-\rho_\T-mu=w\mu-\mu\in\Q$ for all $w\in W(\T)^\prime$.  Since $\Q\cap\aad$ is contained in the root lattice for ${}^{ru}\Sigma_\T$, we have the result.
\,\,\,\,$\Box$

The following result is the main result.
\begin{thm}
We assume that $\T\subsetneq\Pi$ is strictly normal.  Let $\mu\in\aad$ be such that $\rho_\T+\mu$ is dominant integral and regular.  Let $x,y\in W(\T)^\prime$.
Then, we have
\begin{itemize}
\item[(1)] $M_\T(\rho_\T+x\mu)\subseteq M_\T(\rho_\T+y\mu)$ if and only if $y\leqslant_\T x$. 
\item[(2)] If  $y\leqslant_\T x$, then $M_\T(\rho_\T+x\mu)\subseteq M_\T(\rho_\T+y\mu)$ is a composition of  elementary homomorphisms.
\end{itemize}
\end{thm}
\proof 
First, we assume that $M_\T(\rho_\T+x\mu)\subseteq M_\T(\rho_\T+y\mu)$.
Hence, $L(\rho_\T+x\mu)$ is an irreducible constituent of $M(\rho_\T+y\mu)$.
From \cite{[BGG]}, we have$M(\rho_\T+x\mu)\subseteq M(\rho_\T+y\mu)$, namely  $y\leqslant x$.  Hence from Theorem 3.4.5, we have $y\leqslant_\T x$.

Next, we assume that $y\leqslant_\T x$.  Since $\mu$ is regular domimant
integral with respect to ${}^{ru}\Sigma_\T$, there exist  $\alpha_1,...,\alpha_\in{}^{ru}\Delta_\T^+$ such that $\sigma_{\alpha_1}\cdots\sigma_{\alpha_k}y=x$, $\langle y\mu, \alpha_k\rangle >0$, and  $\langle\sigma_{\alpha_{r+1}}\cdots\sigma_{\alpha_k}y\mu, \alpha_r\rangle >0$ for $1\leqslant r\leqslant k-1$.
So, from Proposition 4.2.2, we can construct embedding $M_\T(\rho_\T+x\mu)\subseteq M_\T(\rho_\T+y\mu)$ as a composition of elementary homomorphisms.

\,\,\,\,Q.E.D.
\subsection{Classification of the strictly normal parabolic subalgebras}

From \cite{[M]}, we can determine $\T$-excellent roots and we can obtain the following result.
\begin{prop}
The following is the list of the strictly normal standard parabolic subalgebras of a classical Lie algebra.
\begin{itemize}
\item[(a)]$\ppp(A_{n-1,k})$ \,\,\,\, ($k|n$),
\item[(b)]$\ppp(B_{n,2k,m})$ \,\,\, ($k\leqslant m$),
\item[(c)]$\ppp(B_{n,2k+1,m})$ \,\,\, ($k\geqslant m$),
\item[(d)]$\ppp(C_{n,2k,m})$ \,\,\, ($k\leqslant m$),
\item[(e)]$\ppp(C_{n,2k+1,m})$ \,\,\, ($k\geqslant m$),
\item[(f)]$\ppp(D_{n,2k-1,m})$ \,\,\, ($k\leqslant m$, $2\leqslant m$),
\item[(g)]$\ppp(D_{n,2k,m})$ \,\,\, ($k\geqslant m$,$2\leqslant m$ ),
\item[(h)]$\ppp(D_{n,1,0})$.
\end{itemize}
\end{prop}
Next, we state the classification of strictly normal parabolic subalgebras for exceptional Lie algebras. 
It is obtained by more or less straightforward calculation from  \cite{[M]}. .

\begin{prop}
Let $\ggg$ be an exceptional Lie algebra.  We assume  $\T\subsetneq\Pi$ is normal and $\hbox{{\rm card}} (\Pi-\T)\geqslant 2$.   Moreover, we assume that $\T$ is not strictly normal.  Then $\T$ is  $F_{4,14}$, $E_{7,27}$, or $E_{8, 18}$.
\end{prop}

\remark  If $\hbox{{\rm card}} (\Pi-\T)=1$, then $\pps$ is a maximal parabolic subalgebra.  In this case, the homomorphisms between scalar generalized Verma modules are classified in \cite{[M]}.  So, we neglect them.

\setcounter{section}{6}
\setcounter{subsection}{0}

\section*{\S\,\, 6.\,\,\,\, Normal but not strictly normal case}

\subsection{General results}

We assume that $\ggg$ is simple and $\T\subsetneq\Pi$ is normal but not strictly normal.  We also assume $\pps$ is not a maximal parabolic subalgebra, namely $\hbox{card}(\Pi-\T)\geqslant 2$.
Then, we easily see from the classification that ${}^{ru}\Sigma_\T$ is of the type $\hbox{\rm B}_n$. Here, $n=\dim\aas$.
  Moreover, we assume $\T$ is not of the type $F_{4, 14}$.
Then,  ${}^e\Sigma_\T$ is the set of the long roots in ${}^{ru}\Sigma_\T$.
So, ${}^e\Sigma_\T$ is a root system of the type $\hbox{\rm D}_n$.
We put ${}^{ru}\Pi_\T=\{\gamma_1,\gamma_2,...,\gamma_n\}$ such that $\langle\gamma_i,\gamma_{i+1}^\vee\rangle=-1$ for $1\leqslant i<n$ and $\langle\gamma_2,\gamma_1^\vee\rangle=-2$. Namely, $\gamma_1$ is a unique short simple root.  Put $\gamma^\prime=\sigma_{\gamma_1}\gamma_2$.
Then, $\{\gamma^\prime,\gamma_2,...,\gamma_n\}$ is a basis of ${}^e\Sigma_\T^+$. We put ${}^e S(\T)=\{\sigma_{\gamma^\prime}\}\cup \{\sigma_{\gamma_i}\mid 2\leqslant i\leqslant n\}$ and  we denote by $\leqslant_\T^\prime$ the Bruhat ordering for 
a Coxeter system $({}^e W(\T), {}^e S(\T))$.  We remark that the index of ${}^eW(\T)$ in $W(\T)^\prime$ is two and $W(\T)^\prime={}^eW(\T)\cup {}^eW(\T)\sigma_{\gamma_1}$.

The argument of the proof of Theorem 5.1.5 is partially applicable in this setting and we have the following weaker result.
\begin{prop}
We retain the above setting.
 Let $\mu\in\aad$ be such that $\rho_\T+\mu$ is dominant integral and regular.  
\begin{itemize}
\item[(1)] Let $x,y\in{}^e W(\T)$ be such that $y\leqslant_\T^\prime x$. 
Then, we have $M_\T(\rho_\T+x\mu)\subseteq M_\T(\rho_\T+y\mu)$ and $M_\T(\rho_\T+x\sigma_{\gamma_1}\mu)\subseteq M_\T(\rho_\T+y\sigma_{\gamma_1}\mu)$.  Moreover, $M_\T(\rho_\T+x\mu)\subseteq M_\T(\rho_\T+y\mu)$ and $M_\T(\rho_\T+x\sigma_{\gamma_1}\mu)\subseteq M_\T(\rho_\T+y\sigma_{\gamma_1}\mu)$ are  compositions of  elementary homomorphisms.
\item[(2)] Let $z,w\in W(\T)^\prime$.  If $M_\T(\rho_\T+z\mu)\subseteq M_\T(\rho_\T+w\mu)$, then $w\leqslant_\T z$.
\end{itemize}
\end{prop}

We also have the following result.
\begin{prop}
 Let $\mu\in\aad$ be such that $\rho_\T+\mu$ is dominant integral and regular. 
Let $x,y\in{}^e W(\T)$.
Then, we have $M_\T(\rho_\T+y\sigma_{\gamma_1}\mu)\nsubseteq M_\T(\rho_\T+x\mu)$ and $M_\T(\rho_\T+x\mu)\nsubseteq M_\T(\rho_\T+y\sigma_{\gamma_1}\mu)$.
\end{prop}
\proof
Let $w_1\in{}^eW(\T)$ be the longest element  for  $({}^e W(\T), {}^e S(\T))$.
From proposition 4.4.1, we have $M_\T(\rho_\T+w_1\mu)\subseteq M_\T(\rho_\T+x\mu)$ and  $M_\T(\rho_\T+w_1\sigma_{\gamma_1}\mu)\subseteq M_\T(\rho_\T+y\sigma_{\gamma_1}\mu)$.

We assume that $n$ is even.  Then, $w_1$ is the longest element for $(W(\T)^\prime, S(\T))$.  So, $M_\T(\rho_\T+w_1\mu)$ is irreducible from \cite{[VU]} Proposition 8.5. Since $\gamma_1\not\in{}^e\Sigma_\T$, we have $M_\T(\rho_\T+w_1\mu)\nsubseteq M_\T(\rho_\T+w_1\sigma_{\gamma_1}\mu)$ from Proposition 4.1.4 and $w_1\sigma_{\gamma_1}=\sigma_{\gamma_1}w_1$.
We assume that $M_\T(\rho_\T+w_1\mu)\subseteq M_\T(\rho_\T+y\sigma_{\gamma_1}\mu)$.  The Bernstein degree of a scalar generalized Verma module is one, it contains only one irreducible constituent of the maximal Gelfand-Kirillov dimension.  So, from Proposition 1.4.1  $M_\T(\rho_\T+w_1\sigma_{\gamma_1}\mu)\subseteq M_\T(\rho_\T+y\sigma_{\gamma_1}\mu)$ implies that $M_\T(\rho_\T+w_1\mu)\subseteq M_\T(\rho_\T+w_1\sigma_{\gamma_1}\mu)$.
So, we have a contradiction.  Hence, $M_\T(\rho_\T+w_1\mu)\nsubseteq M_\T(\rho_\T+y\sigma_{\gamma_1}\mu)$. Since $M_\T(\rho_\T+w_1\mu)\subseteq M_\T(\rho_\T+x\mu)$, we have $M_\T(\rho_\T+x\mu)\nsubseteq M_\T(\rho_\T+y\sigma_{\gamma_1}\mu)$.
Next, we assume that $M_\T(\rho_\T+y\sigma_{\gamma_1}\mu)\subseteq M_\T(\rho_\T+x\mu)$.   $M_\T(\rho_\T+w_1\mu)$ is the unique irreducible constituent of $M_\T(\rho_\T+x\mu)$ of the maximal Gelfand-Kirillov dimension from Proposition 1.4.1 and $M_\T(\rho_\T+w_1\mu)\nsubseteq M_\T(\rho_\T+y\sigma_{\gamma_1}\mu)$.
So, we have $M_\T(\rho_\T+w_1\mu)\subseteq M_\T(\rho_\T+x\mu)/M_\T(\rho_\T+y\sigma_{\gamma_1}\mu)$.
On the other hand, $\Dim(M_\T(\rho_\T+x\mu)/M_\T(\rho_\T+y\sigma_{\gamma_1}\mu))<\Dim(M_\T(\rho_\T+w_1\mu))$
from Proposition 1.4.1.
So, we have a contradiction.
This means that $M_\T(\rho_\T+y\sigma_{\gamma_1}\mu)\nsubseteq M_\T(\rho_\T+x\mu)$. 

If we $n$ is odd,then  $w_1\sigma_{\gamma_1}$ is the longest element for $(W(\T)^\prime, S(\T))$. The proof of the proposition in this case is more or less similar to that for the case that $n$ is even.  So, we omit proving the proposition in this case.
 \,\,\,\,$\Box$

\subsection{$B_2$ case}
Next, we consider the case of $n=2$.  We fix $\mu\in\aad$ such that $\rho_\T+\mu$ is regular dominant integral.
In this case, ${}^eW(\T)=\{e, \sigma_{\gamma_2},  \sigma_{\gamma^\prime}, \sigma_{\gamma_2} \sigma_{\gamma^\prime}\}\cong \itg/2\itg\times\itg/2\itg$.
Since $\sigma_{\gamma_2}\leqslant_\T \sigma_{\gamma^\prime}$ and $\sigma_{\gamma_2}\nleqslant_\T ^\prime\sigma_{\gamma^\prime}$ hold, Proposition 6.1.1 and Proposition 6.1.2 are insufficient to determine whether 
$M_\T(\rho_\T+\sigma_{\gamma^\prime}\mu)\subseteq M_\T(\rho_\T+ \sigma_{\gamma_2}\mu)$ or not.
However, as a corollary of Proposition 6.1.1 and Proposition 6.1.2, we easily have the following result.
\begin{cor}
Let $x,y\in W(\T)^\prime$ such that $x\neq y$ and $(x,y)\neq (\sigma_{\gamma^\prime}, \sigma_{\gamma_2})$.
Then, we have
\begin{itemize}
\item[(1)] $M_\T(\rho_\T+x\mu)\subseteq M_\T(\rho_\T+y\mu)$ if and only if $(x,y)$ appears in the following list.

 $(\sigma_{\gamma_2},e), (\sigma_{\gamma^\prime},e), (\sigma_{\gamma_2} \sigma_{\gamma^\prime},e), (\sigma_{\gamma_2} \sigma_{\gamma^\prime},
\sigma_{\gamma^\prime}), (\sigma_{\gamma_2} \sigma_{\gamma^\prime}, \sigma_{\gamma_2}), 
 (\sigma_{\gamma_2}\sigma_{\gamma_1},\sigma_{\gamma_1}),$

$ (\sigma_{\gamma^\prime}\sigma_{\gamma_1},\sigma_{\gamma_1}), (\sigma_{\gamma_2} \sigma_{\gamma^\prime}\sigma_{\gamma_1},\sigma_{\gamma_1}), (\sigma_{\gamma_2} \sigma_{\gamma^\prime}\sigma_{\gamma_1},\sigma_{\gamma^\prime}\sigma_{\gamma_1}), (\sigma_{\gamma_2} \sigma_{\gamma^\prime}\sigma_{\gamma_1}, \sigma_{\gamma_2}\sigma_{\gamma_1}).
$
\item[(2)] If $M_\T(\rho_\T+x\mu)\subseteq M_\T(\rho_\T+y\mu)$, then it is a composition of elementary homomorphisms.
\end{itemize}
\end{cor}

\subsection{$\hbox{\rm B}_{n,1,n-2}$}
We fix notation as follows. 
Let $n\geqslant 3$.
Let $\gggg$ be a simple Lie algebra of the type $\hbox{\rm B}_{n}$.
We can choose an orthonormal basis $e_1,...,e_n$ of $\hhh^\ast$ such that 
\begin{align*}
\Delta = \{ \pm e_i\pm e_j\mid 1\leqslant i< j\leqslant
n\}\cup
\{\pm e_i\mid 1\leqslant i\leqslant n\}.
\end{align*}
We choose a positive system as follows.
\begin{align*}
\Delta^+ = \{ e_i\pm e_j\mid 1\leqslant i< j\leqslant
n\}\cup
\{e_i\mid 1\leqslant i\leqslant n\}.
\end{align*}
If we put $\alpha_i=e_i-e_{i+1}$  \,\,\, $(1\leqslant i <n)$ and
$\alpha_n=e_n$,
then $\Pi=\{\alpha_1,...,\alpha_n\}$.
We put $\T=\Pi-\{\alpha_1, \alpha_2\}$.
$$
\bigcirc
-\bigcirc
-\textcircled{$\bullet$}-\textcircled{$\bullet$}-\cdots-\textcircled{$\bullet$}\Rightarrow\textcircled{$\bullet$}
$$
Then, we have 
$\aad=\{se_1+te_2\mid s,t\in\cpx\}$ and $\rho_\T=\sum_{i=3}^{m}\frac{2m-2i+1}{2}e_{i}$.
We put $\gamma_1=\alpha_2|_{\aas}=e_2$, $\gamma_2=\alpha_1|_{\aas}=e_1-e_2$, and $\gamma^\prime=\sigma_{\gamma_1}\gamma_2=
e_1+e_2$.
Then, these notations are compatible with those in 6.2 and 6.3.

we put $\nu_0=\frac{3}{2}e_1+\frac{1}{2}e_2$
Then, $\rho_\T+\nu_0$ is integral and $\langle \nu_0,\gamma\rangle\geqslant 0$ for all $\gamma\in\Sigma_\T^+$.
We put $\mu_1=\frac{1}{2}e_1+\frac{1}{2}e_2$, 
$\mu_2=\frac{1}{2}e_1-\frac{1}{2}e_2$,
$\mu_3=-\frac{1}{2}e_1+\frac{1}{2}e_2$, and
$\mu_4=-\frac{1}{2}e_1-\frac{1}{2}e_2$.

First, we have the following results by a straightforward computation.
(Note: $\rho_\T+\sigma_{\gamma_2}\nu_0=\frac{1}{2}e_1+\frac{3}{2}e_2+\sum_{i=3}^{m}\frac{2m-2i+1}{2}e_{i}$, $\rho_\T+\sigma_{\gamma_2}\sigma_{\gamma_1}\nu_0=-\frac{1}{2}e_1+\frac{3}{2}e_2+\sum_{i=3}^{m}\frac{2m-2i+1}{2}e_{i}=\rho_\T+\sigma_{\gamma_2}\nu_0-e_1$.)
\begin{lem}
$$\{\rho_\T+\sigma_{\gamma_2}\nu_0\pm e_i\mid 1\leqslant k\leqslant n\}\cap\PS\cap W\dot(\rho_\T+\sigma_{\gamma_2}\nu_0)
=\{ \rho_\T+\sigma_{\gamma_2}\sigma_{\gamma_1}\nu_0\}.$$
\end{lem}

Let $V$ be a natural representation of $\ggg$.  Namely, $V$ is an irreducible representation of $V$ with a highest weight $e_1$.
	The, the set of the weights of $V$ is $\{\pm e_i\mid 1\leqslant k\leqslant n\}\cup\{0\}$.
We also easily have the following result.
\begin{lem}

There exists a sequence of $\pps$-submodules of $V$ 
$$ \{0\}=V_0\subsetneq V_1\subsetneq V_2\subsetneq\cdots\subsetneq V_{2n+1}=V$$
which satisfies the following conditions (a)-(e).
\begin{itemize}
\item[(a)] $V_{i}/V_{i-1}$ is a one-dimensional $\pps$-module such that $\nns$ acts on it trivially for each $1\leqslant i\leqslant 2n+1$.
\item[(b)] We denote by $\lambda_i$ the highest weight of $V_{i}/V_{i-1}$ as an $\lls$-module.  Then $\lambda_i=e_i$ for $1\leqslant i\leqslant n$, $\lambda_{n+1}=0$, and $\lambda_i=-e_{2n+2-i}$ for $n+2\leqslant i\leqslant 2n+1$.
\end{itemize}
\end{lem}

\begin{lem}
\begin{itemize}
\item[(1)] $M_\T(\rho_\T+\mu_3)$ is irreducible.
\item[(2)]$M_\T(\rho_\T+\mu_4)$ is irreducible.
\end{itemize}
\end{lem}
\proof
(1) can be  proved by Janzten's irreducibility condition (\cite{[Ja]} Satz 3).
(2) follows from \cite{[VU]} Proposition 8.5.
\,\,\,\, $\Box$
\begin{lem}
$$M_\T(\rho_\T+\sigma_{\gamma^\prime}\sigma_{\gamma_1}\nu_0)\nsubseteq
M_\T(\rho_\T+ \sigma_{\gamma_2}\nu_0).$$
\end{lem}
\proof
We assume that $M_\T(\rho_\T+\sigma_{\gamma^\prime}\sigma_{\gamma_1}\nu_0)\subseteq
M_\T(\rho_\T+ \sigma_{\gamma_2}\nu_0)$ and deduce a contradiction.

From Lemma 2.4.6 and the exactness of a translation functor, we have 
$$M_\T(\rho_\T+\mu_2)\subseteq
M_\T(\rho_\T+\mu_1).$$
From Proposition 4.2.2, we have $M_\T(\rho_\T+\mu_3)\subseteq 
M_\T(\rho_\T+\mu_2)$.  So, we have $M_\T(\rho_\T+\mu_3)\subseteq M_\T(\rho_\T+\mu_1)$.
From Proposition 4.2.2, we have $M_\T(\rho_\T+\mu_4)\subseteq M_\T(\rho_\T+\mu_1)$.
Hence $M_\T(\rho_\T+\mu_1)$ has at least two distinct irreducible submodules.  It contradicts Proposition 1.4.1 (3).  \,\,\,\, $\Box$

\begin{lem}
$$ M_\T(\rho_\T+\sigma_{\gamma^\prime}\nu_0)\nsubseteq M_\T(\rho_\T+ \sigma_{\gamma_2}\nu_0).$$
\end{lem}

\proof
Assuming $M_\T(\rho_\T+\sigma_{\gamma^\prime}\nu_0)\subseteq M_\T(\rho_\T+ \sigma_{\gamma_2}\nu_0)$, we deduce a contradiction.
Put $X=P_{\rho_\T+\nu_0}(M_\T(\rho_\T+\sigma_{\gamma^\prime}\nu_0)\otimes V)$ and  $Y=P_{\rho_\T+\nu_0}(M_\T(\rho_\T+\sigma_{\gamma_2}\nu_0)\otimes V)$.
Hence, we have $X\subseteq Y$.

We remark that $\sigma_{\gamma_2}\sigma_{\gamma_1}\nu_0-\sigma_{\gamma_2}\nu_0=-e_1$. 
So, from Lemma 6.3.1 and Lemma 6.3.2, we see that  there is a  submodule $Y_1\subsetneq Y$ such that $Y_1\cong M_\T(\rho_\T+\sigma_{\gamma_2}\nu_0)$ and $Y/Y_1\cong M_\T(\rho_\T+\sigma_{\gamma_2}\sigma_{\gamma_1}\nu_0)$.
On the other hand,  $\sigma_{\gamma^\prime}\sigma_{\gamma_1}\nu_0-\sigma_{\gamma^\prime}\nu_0=e_1$.
and it is the highest weight of the natural representation $V$.
So, there is an embedding 
$$ \iota: M_\T(\rho_\T+\sigma_{\gamma^\prime}\sigma_{\gamma_1}\nu_0)\subseteq X\subseteq Y.$$
Considering composition with the canonical projection
$$q: Y\rightarrow Y/Y_1\cong M_\T(\rho_\T+\sigma_{\gamma_2}\sigma_{\gamma_1}\nu_0),$$
we obtain a homomorphism
 $$q\circ\iota: M_\T(\rho_\T+\sigma_{\gamma^\prime}\sigma_{\gamma_1}\nu_0)\rightarrow
M_\T(\rho_\T+\sigma_{\gamma_2}\sigma_{\gamma_1}\nu_0).$$
We assume that $q\circ\iota=0$.  Then, we have $M_\T(\rho_\T+\sigma_{\gamma^\prime}\sigma_{\gamma_1}\nu_0)\subseteq Y_1$. 
However, Lemma 6.3.4 implies that it is a zero map.
So, we have $q\circ\iota\neq 0$.  From Theorem 2.1.1 (2), we have $M_\T(\rho_\T+\sigma_{\gamma^\prime}\sigma_{\gamma_1}\nu_0)\subseteq
M_\T(\rho_\T+\sigma_{\gamma_2}\sigma_{\gamma_1}\nu_0).$
From Proposition 2.4.6 and the exactness of the translation functors, we have 
$M_\T(\rho_\T+\mu_2)\subseteq 
M_\T(\rho_\T+\mu_3)$.  However, in our proof of Lemma 6.3.4, we see $M_\T(\rho_\T+\mu_3)\subseteq 
M_\T(\rho_\T+\mu_2)$.
So, we obtained a contradiction.  \,\,\,\, $\Box$

\begin{prop}
Let $\mu\in\aad$ be such that $\rho_\T+\mu$ is dominant integral and regular.  
$$ M_\T(\rho_\T+\sigma_{\gamma^\prime}\mu)\nsubseteq M_\T(\rho_\T+ \sigma_{\gamma_2}\mu).$$
\end{prop}

\proof We assume that $ M_\T(\rho_\T+\sigma_{\gamma^\prime}\mu)\subseteq M_\T(\rho_\T+ \sigma_{\gamma_2}\mu).$  From the translation principle, we have 
$ M_\T(\rho_\T+\sigma_{\gamma^\prime}\rho^\T)\subseteq M_\T(\rho_\T+ \sigma_{\gamma_2}\rho^\T).$
Since $\rho^\T-\nu_0$ is dominant integral, we have 
$ M_\T(\rho_\T+\sigma_{\gamma^\prime}\nu_0)\subseteq M_\T(\rho_\T+ \sigma_{\gamma_2}\nu_0)$  from  Lemma 2.4.7.  It contradicts Lemma 6.3.5.  \,\,\,\, $\Box$

\begin{cor}
Let $\lambda,\mu\in\aad$ be such that $\rho_\T+\mu$ and $\rho_\T+\lambda$ is integral and regular and $M_\T(\rho_\T+\mu)\subseteq M_\T(\rho_\T+ \lambda)$. 
Then,  $M_\T(\rho_\T+\mu)\subseteq M_\T(\rho_\T+ \lambda)$ is a composition of some elementary homomorphisms.
\end{cor}

\subsection{$\hbox{\rm E}_{7,27}$ and $\hbox{\rm E}_{8,18}$}

We fix notation as follows. 

First, we consider $\hbox{\rm E}_{7,27}$.  Let $\gggg$ be a simple Lie algebra of the type $E_7$.

We fix an orthonormal basis $e_1,...,e_8$ in $\rel^8$.
We identify $\hhh^\ast$ with $\{v\in\rel^8\mid \langle v,
e_1-e_2\rangle=0\}$ so that
\begin{multline*}
\Delta=\{\pm (e_1+e_2)\}\cup\{\pm e_i\pm e_j\mid 3\leqslant i<j\leqslant
8\}\\
\cup \left\{\pm\frac{1}{2}\left(\left. e_1+e_2+\sum_{i=3}^8\varepsilon_i e_i\right)\right|
\mbox{$\varepsilon_i=\pm 1$ is for $3\leqslant i\leqslant 8$ and
 $\prod_{i=3}^8\varepsilon_i=1$.} \right\}
\end{multline*}

We choose a positive system as follows.
\begin{multline*}
\Delta^+=\{ (e_1+e_2)\}\cup\{ e_i\pm e_j\mid 3\leqslant i<j\leqslant
8\}\\
\cup \left\{\frac{1}{2}\left(\left. e_1+e_2+\sum_{i=3}^8\varepsilon_i e_i\right)\right|
\mbox{$\varepsilon_i=\pm 1$ for $3\leqslant i\leqslant 8$ and
 $\prod_{i=3}^8\varepsilon_i=1$.} \right\}
\end{multline*}

Put $\alpha_i=e_{i+2}-e_{i+3}$ for $1\leqslant i\leqslant 5$,
$\alpha_6=e_7+e_8$, and
$\alpha_7=\frac{1}{2}(e_1+e_2-e_3-e_4-e_5-e_6-e_7-e_8)$.
Then, $\Pi=\{\alpha_1,...,\alpha_7\}$ is the set of simple roots in
$\Delta^+$.

\[
 \begin{array}{ccccccccccc} 1 & - & 2 & - & 3 & - & 4 & - & 6 & - & 7\\
& &  & & & & \mid & & & & \\
& &  & & & & 5 & & & & 
\end{array}
\]
We consider the standard parabolic subalgebra of the type $\hbox{\rm E}_{7,27}$.
Namely, we put $\T=\Pi-\{\alpha_2,\alpha_7\}$.
Then, we have 
$\aad=\{se_1+se_2+te_3+te_4\mid s,t\in\cpx\}$ and $\rho_\T=\frac{1}{2}e_3-\frac{1}{2}e_4+3e_5+2e_6+e_7$.
We put $\gamma_1=\alpha_2|_{\aas}=\frac{1}{2}e_3+\frac{1}{2}e_4$, $\gamma_2=\alpha_7|_{\aas}=\frac{1}{2}e_1+\frac{1}{2}e_2-\frac{1}{2}e_3-\frac{1}{2}e_4$, and $\gamma^\prime=\sigma_{\gamma_1}\gamma_2=
\frac{1}{2}e_1+\frac{1}{2}e_2+\frac{1}{2}e_3+\frac{1}{2}e_4$.
Then, these notations are compatible with those in 6.2 and 6.3.
we put $\nu_0=\frac{3}{2}e_1+\frac{3}{2}e_2-\frac{1}{2}e_3-\frac{1}{2}e_4$
Then, $\rho_\T+\nu_0=\frac{3}{2}e_1+\frac{3}{2}e_2-e_4+3e_5+2e_6+e_7$ is integral and $\rho-(\rho_\T+\nu_0)=7e_1+7e_2+5e_3+3e_4+3e_5+2e_6+e_7$ is dominant integral.
We put $\mu_1=\frac{1}{2}e_1+\frac{1}{2}e_2+\frac{1}{2}e_3+\frac{1}{2}e_4$, 
$\mu_2=\frac{1}{2}e_1+\frac{1}{2}e_2-\frac{1}{2}e_3-\frac{1}{2}e_4$,
$\mu_3=-\frac{1}{2}e_1-\frac{1}{2}e_2+\frac{1}{2}e_3+\frac{1}{2}e_4$,
$\mu_4=-\frac{1}{2}e_1-\frac{1}{2}e_2-\frac{1}{2}e_3-\frac{1}{2}e_4$.

Next, we consider the case of $\hbox{\rm E}_{8,18}$.
Let $\gggg$ be a simple Lie algebra of the type $E_7$.
We fix an orthonormal basis $e_1,...,e_8$ in $\hhh^\ast$
such  that
\begin{multline*}
\Delta=\{\pm e_i\pm e_j\mid 1\leqslant i<j\leqslant
8\}\\
\cup \left\{\pm\frac{1}{2}\left(\left.\sum_{i=1}^8\varepsilon_i e_i\right)\right|
\mbox{$\varepsilon_i=\pm 1$for $1\leqslant i\leqslant 8$ and
 $\prod_{i=1}^8\varepsilon_i=-1$.} \right\}
\end{multline*}

We choose a positive system as follows.
\begin{multline*}
\Delta^+=\{ e_i\pm e_j\mid 1\leqslant i<j\leqslant
8\}\\
\cup \left\{\frac{1}{2}\left(\left. e_1+\sum_{i=2}^8\varepsilon_i e_i\right)\right|
\mbox{$\varepsilon_i=\pm 1$ is for $2\leqslant i\leqslant 8$ and
 $\prod_{i=2}^8\varepsilon_i=-1$.} \right\}
\end{multline*}

Put $\alpha_i=e_{i+1}-e_{i+2}$ for $1\leqslant i\leqslant 6$,
$\alpha_7=e_7+e_8$, and
$\alpha_8=\frac{1}{2}(e_1-e_2-e_3-e_4-e_5-e_6-e_7-e_8)$.
Then, $\Pi=\{\alpha_1,...,\alpha_8\}$ is the set of simple roots in
$\Delta^+$.
\[
 \begin{array}{ccccccccccccc} 1 & - & 2 & - & 3 & - & 4 & - & 5 & - & 7 & - & 8\\
& &  & & & & & & \mid & & & & \\
& & & &  & & & & 6 & & & & 
\end{array}
\]

We consider the standard parabolic subalgebra of the type $\hbox{\rm E}_{8,18}$.
Namely, we put $\T=\Pi-\{\alpha_1,\alpha_8\}$.
Then, we have 
$\aad=\{se_1+te_2\mid s,t\in\cpx\}$ and $\rho_\T=5e_3+4e_4+3e_5+2e_6+e_7$.
We put $\gamma_1=\alpha_8|_{\aas}=e_1-e_2$, $\gamma_2=\alpha_1|_{\aas}=e_2$, and $\gamma^\prime=\sigma_{\gamma_1}\gamma_2=
e_1$.
Then, these notations are compatible with those in 6.2 and 6.3.
we put $\nu_0=2e_1+e_2$.  It is dominant with respect to $\Sigma_\T^+$.
Then, $\rho_\T+\nu_0=2e_1+e_2+5e_3+4e_4+3e_5+2e_6+e_7$ is integral and $\rho-(\rho_\T+\nu_0)=21e_1+5e_2$ is dominant integral.
We put $\mu_1=e_1$, 
$\mu_2=-e_2$,
$\mu_3=e_2$,
$\mu_4=-e_1$.

Hereafter, we assume that $\ggg$ is  of either the type $\hbox{\rm E}_{7}$ or $\hbox{\rm E}_{8}$.  We treat these case at the same time.

First, we have the following results by a straightforward computation.
\begin{lem}
\begin{multline*}
\{\rho_\T+\sigma_{\gamma_2}\nu_0+\beta\mid \beta\in\Delta\cup\{0\}\}\cap\PS\cap W\dot(\rho_\T+\sigma_{\gamma_2}\nu_0)\\
=\{\rho_\T+\sigma_{\gamma_2}\nu_0, \rho_\T+\sigma_{\gamma_2}\sigma_{\gamma_1}\nu_0\}.
\end{multline*}
\end{lem}
We also have the following result.
\begin{lem}
We regard $\ggg$ as an adjoint representation.

\begin{itemize}
\item[(1)]
There exists a sequence of $\pps$-submodules of $\ggg$ 
$$ \{0\}=V_0\subsetneq V_1\subsetneq V_2\subsetneq\cdots\subsetneq V_k=\ggg$$
and a sequence of roots  $\beta_1,...,\beta_k\in\Delta$
which satisfies the following conditions (a)-(e).
\begin{itemize}
\item[(a)] $V_{i}/V_{i-1}$ is an irreducible $\pps$-module such that $\nns$ acts on it trivially for each $1\leqslant i\leqslant k$.
\item[(b)] As a $\lls$-module, $V_{i}/V_{i-1}$ has a highest weight $\beta_i$ for each $1\leqslant i\leqslant k$.
\item[(c)] $\beta_1$ is the highest root $e_1+e_2$.  
\item[(d)] $\beta_k$ is the lowest root $-e_1-e_2$.
\item[(e)] There exist some $1<h_1<h_2<k$ such that $\beta_{h_1}=\beta_{h_2}=0$ and $\beta_i\neq 0$ for all $i\neq h_1.h_2$.
\end{itemize}
\item[(2)] $V_1$, $V_{h_1}/V_{h_1-1}$, $V_{h_2}/V_{h_2-1}$, $V_{k}/V_{k-1}$
are all one-dimensional.
\end{itemize}
\end{lem}
\proof
The existence of $V_1,...,V_k$ satisfying (a)-(d) is proved by a standard argument.
For $H\in\hhh$, $[\lls\cap\nnn, H]=0$ if and  only if $H\in\aas$.
Since $\dim\aas=2$, $V_1,...,V_k$ satisfies (e) also.
From $\pm(e_1+e_2), 0\in\aad$, we easily have (2).
\,\,\,\, $\Box$

\begin{lem}
\begin{itemize}
\item[(1)] $M_\T(\rho_\T+\mu_3)$ is irreducible.
\item[(2)]$M_\T(\rho_\T+\mu_4)$ is irreducible.
\end{itemize}
\end{lem}
\proof
(1) can be  proved by Janzten's irreducibility condition (\cite{[Ja]} Satz 3).
(2) follows from \cite{[VU]} Proposition 8.5.
\,\,\,\, $\Box$
\begin{lem}
$$M_\T(\rho_\T+\sigma_{\gamma^\prime}\sigma_{\gamma_1}\nu_0)\nsubseteq
M_\T(\rho_\T+ \sigma_{\gamma_2}\nu_0).$$
\end{lem}
\proof
We assume that $M_\T(\rho_\T+\sigma_{\gamma^\prime}\sigma_{\gamma_1}\nu_0)\subseteq
M_\T(\rho_\T+ \sigma_{\gamma_2}\nu_0)$ and deduce a contradiction.

From Lemma 2.4.6 and the exactness of a translation functor, we have 
$$M_\T(\rho_\T+\mu_2)\subseteq
M_\T(\rho_\T+\mu_1).$$
From Proposition 4.2.2, we have $M_\T(\rho_\T+\mu_3)\subseteq 
M_\T(\rho_\T+\mu_2)$.  So, we have $M_\T(\rho_\T+\mu_3)\subseteq M_\T(\rho_\T+\mu_1)$.
From Proposition 4.2.2, we have $M_\T(\rho_\T+\mu_4)\subseteq M_\T(\rho_\T+\mu_1)$.
Hence $M_\T(\rho_\T+\mu_1)$ has at least two distinct irreducible submodules.  It contradicts Proposition 1.4.1 (3).  \,\,\,\, $\Box$

\begin{lem}
$$ M_\T(\rho_\T+\sigma_{\gamma^\prime}\nu_0)\nsubseteq M_\T(\rho_\T+ \sigma_{\gamma_2}\nu_0).$$
\end{lem}

\proof
Assuming $M_\T(\rho_\T+\sigma_{\gamma^\prime}\nu_0)\subseteq M_\T(\rho_\T+ \sigma_{\gamma_2}\nu_0)$, we deduce a contradiction.
Put $X=P_{\rho_\T+\nu_0}(M_\T(\rho_\T+\sigma_{\gamma^\prime}\nu_0)\otimes\ggg)$ and  $Y=P_{\rho_\T+\nu_0}(M_\T(\rho_\T+\sigma_{\gamma_2}\nu_0)\otimes\ggg)$.
Hence, we have $X\subseteq Y$.

We remark that $\sigma_{\gamma_2}\sigma_{\gamma_1}\nu_0-\sigma_{\gamma_2}\nu_0=-e_1-e_2$. 
So, from Lemma 6.4.1 and Lemma 6.4.2, we see that  there is a sequence of submodules $Y_1\subsetneq Y_2\subsetneq Y$ such that $Y_1\cong Y_2/Y_1\cong M_\T(\rho_\T+\sigma_{\gamma_2}\nu_0)$ and $Y/Y_2\cong M_\T(\rho_\T+\sigma_{\gamma_2}\sigma_{\gamma_1}\nu_0)$.
On the other hand,  $\sigma_{\gamma^\prime}\sigma_{\gamma_1}\nu_0-\sigma_{\gamma^\prime}\nu_0=e_1+e_2$.
and it is the highest weight of the adjoint representation $\ggg$.
So, there is an embedding 
$$ \iota: M_\T(\rho_\T+\sigma_{\gamma^\prime}\sigma_{\gamma_1}\nu_0)\subseteq X\subseteq Y.$$
Considering composition with the canonical projection
$$q: Y\rightarrow Y/Y_2\cong M_\T(\rho_\T+\sigma_{\gamma_2}\sigma_{\gamma_1}\nu_0),$$
we obtain a homomorphism
 $$q\circ\iota: M_\T(\rho_\T+\sigma_{\gamma^\prime}\sigma_{\gamma_1}\nu_0)\rightarrow
M_\T(\rho_\T+\sigma_{\gamma_2}\sigma_{\gamma_1}\nu_0).$$
We assume that $q\circ\iota=0$.  Then, we have $M_\T(\rho_\T+\sigma_{\gamma^\prime}\sigma_{\gamma_1}\nu_0)\subseteq Y_2$.  Again we consider the composition $M_\T(\rho_\T+\sigma_{\gamma^\prime}\sigma_{\gamma_1}\nu_0)\subseteq Y_2\rightarrow Y_2/Y_1\cong  M_\T(\rho_\T+\sigma_{\gamma_2}\nu_0)$.
However, Lemma 6.4.4 implies that it is a zero map.
So, we have $M_\T(\rho_\T+\sigma_{\gamma^\prime}\sigma_{\gamma_1}\nu_0)\subseteq Y_1
\cong  M_\T(\rho_\T+\sigma_{\gamma_2}\nu_0)$.  It contradicts Lemma 6.3.4,
So, we have $q\circ\iota\neq 0$.  From Theorem 2.1.1 (2), we have $M_\T(\rho_\T+\sigma_{\gamma^\prime}\sigma_{\gamma_1}\nu_0)\subseteq
M_\T(\rho_\T+\sigma_{\gamma_2}\sigma_{\gamma_1}\nu_0).$
From Proposition 2.4.6 and the exactness of the translation functors, we have 
$M_\T(\rho_\T+\mu_2)\subseteq 
M_\T(\rho_\T+\mu_3)$.  However, in our proof of Lemma 6.4.4, we see $M_\T(\rho_\T+\mu_3)\subseteq 
M_\T(\rho_\T+\mu_2)$.
So, we obtained a contradiction.  \,\,\,\, $\Box$

\begin{prop}
Let $\mu\in\aad$ be such that $\rho_\T+\mu$ is dominant integral and regular.  
$$ M_\T(\rho_\T+\sigma_{\gamma^\prime}\mu)\nsubseteq M_\T(\rho_\T+ \sigma_{\gamma_2}\mu).$$
\end{prop}

\proof We assume that $ M_\T(\rho_\T+\sigma_{\gamma^\prime}\mu)\subseteq M_\T(\rho_\T+ \sigma_{\gamma_2}\mu).$  From the translation principle, we have 
$ M_\T(\rho_\T+\sigma_{\gamma^\prime}\rho^\T)\subseteq M_\T(\rho_\T+ \sigma_{\gamma_2}\rho^\T).$
Since $\rho^\T-\nu_0$ is dominant integral, we have 
$ M_\T(\rho_\T+\sigma_{\gamma^\prime}\nu_0)\subseteq M_\T(\rho_\T+ \sigma_{\gamma_2}\nu_0)$  from  Lemma 2.4.7.  It contradicts Lemma 6.4.5.  \,\,\,\, $\Box$

\begin{cor}
Let $\lambda,\mu\in\aad$ be such that $\rho_\T+\mu$ and $\rho_\T+\lambda$ is integral and regular and $M_\T(\rho_\T+\mu)\subseteq M_\T(\rho_\T+ \lambda)$. 
Then,  $M_\T(\rho_\T+\mu)\subseteq M_\T(\rho_\T+ \lambda)$ is a composition of some elementary homomorphisms.
\end{cor}

\subsection{$\hbox{\rm F}_{4,14}$}

We consider the root system $\Delta$ for a simple Lie
algebra $\gggg$ of the type $F_4$. (For example, see \cite{[K2]} p691.)
We can choose an orthonormal basis $e_1,...,e_4$ of $\hhh^\ast$ such that 
\begin{align*}
\Delta = \{ \pm e_i\pm e_j\mid 1\leqslant i< j\leqslant
4\}\cup
\{\pm e_i\mid 1\leqslant i\leqslant 4\}
\cup\left\{\frac{1}{2}(\pm e_1\pm e_2\pm e_3 \pm e_4 )\right\}.
\end{align*}
We choose a positive system as follows.
\begin{align*}
\Delta^+ = \{ e_i\pm e_j\mid 1\leqslant i< j\leqslant
4\}\cup
\{e_i\mid 1\leqslant i\leqslant 4\}\cup
\left\{\frac{1}{2}(e_1\pm e_2\pm e_3 \pm e_4 )\right\}.
\end{align*}
Put $\alpha_1=\frac{1}{2}(e_1-e_2-e_3-e_4 )$, $\alpha_2=e_4$,
$\alpha_3=e_3-e_4$, and $\alpha_4=e_2-e_3$.
Then, $\Pi=\{\alpha_1,...,\alpha_4\}$.
\[
 \begin{array}{ccccccc} 1 & - & 2 & \Leftarrow & 3 & - & 4
\end{array}
\]

We consider the standard parabolic subalgebra of the type $\hbox{\rm F}_{4,14}$.
Namely, we put $\T=\Pi-\{\alpha_1,\alpha_4\}$.
Then, we have 
$\aad=\{se_1+te_2\mid s,t\in\cpx\}$ and $\rho_\T=\frac{3}{2}e_3-\frac{1}{2}e_4$.
We put $\gamma_1=\alpha_1|_{\aas}=\frac{1}{2}e_1-\frac{1}{2}e_2$, $\gamma_2=\alpha_4|_{\aas}=e_2$.
In this setting, ${}^{ur}\Sigma_\T$ is of the type $\hbox{\rm B}_2$.
On  ${}^e\Sigma_\T=\emptyset$ and ${}^eW(\T)=\{e\}$. 
We put $\mu_1=\frac{1}{2}e_1+\frac{1}{2}e_2$, 
$\mu_2=\frac{1}{2}e_1-\frac{1}{2}e_2$,
$\mu_3=-\frac{1}{2}e_1+\frac{1}{2}e_2$,
$\mu_4=-\frac{1}{2}e_1-\frac{1}{2}e_2$.
From Janzten's irreducibility criterion, we can show the following result.
\begin{lem}
$M_\T(\rho_\T+\mu_i)$ is irreducible for each $1\leqslant i\leqslant 4$.
\end{lem}
Finally, we prove the following result.

\begin{prop}
Let $\lambda,\nu\in\aad$ be such that $\rho_\T+\lambda$ and $\rho_\T+\nu$ are regular integral and $\lambda\neq\nu$.
Then $M_\T(\rho_\T+\nu)\nsubseteq M_\T(\rho_\T+\lambda)$.
\end{prop}

\proof
We assume that $M_\T(\rho_\T+\nu)\subseteq M_\T(\rho_\T+\lambda)$.
From Proposition 2.3.2, there exists some $x\in W(\T)$ such that $\nu=x\lambda$.
We easy to see there is unique $1\leqslant i\leqslant 4$ (resp.\ $1\leqslant i\leqslant 4$) such that $\lambda$ and $\mu_i$ (resp.\ $x\lambda$ and $\mu_j$) satisfy the condition (T) in 2.4.  From proposition 2.4.7, we have $M_\T(\rho_\T+\mu_j)\subseteq M_\T(\rho_\T+\mu_i)$.  From Lemma 6.5.1, we have $i=j$.  Then, we may apply Proposition 4.1.4 and $M_\T(\rho_\T+\nu)\subseteq M_\T(\rho_\T+\lambda)$is an elementary homomorphism.  However, it contradicts ${}^e\Sigma_\T=\emptyset$.
\,\,\,\,$\Box$

\setcounter{section}{7}
\setcounter{subsection}{0}
\section*{\S\,\, 7.\,\,\,\, Type A case}

\subsection{Some notations}
In this section, we assume that $\ggg=\ggg\llll(n,\cpx)$.  Let $\hhh$ be the Cartan subalgebra of $\ggg$ consisting of the diagonal matrices and let $\bbb$ be the Borel subalgebra of $\ggg$ consisting of the upper triangular matrices. We choose $\Delta^+$ corresponding to $\bbb$.
Then we can choose an orthonormal basis $e_1,...,e_n$ of $\hhh^\ast$ such that 
\begin{align*}
\Delta^+ = \{ e_i-e_j\mid 1\leqslant i< j\leqslant
n, i\neq j\}.
\end{align*}

If we put $\alpha_i=e_i-e_{i+1}$  \,\,\, $(1\leqslant i <n)$ ,
then $\Pi=\{\alpha_1,...,\alpha_{n-1}\}$.
We identify the Weyl group $W$ with the $n$-th symmetric group ${\mathfrak S}_n$ via $\sigma e_i=e_{\sigma(i)}$  \,\,\,\,($1\leqslant i\leqslant n$).

\subsection{Almost normal parabolic subalgebras}

We fix $\T\subsetneq\Pi$.
Put $\T_0=\{\beta\in\T\mid w\beta=\beta \,\, (w\in W(\T))\}$ and put $\T_1=\T-\T_0$. 
We put $\Phi_{\T_0}=\Delta\cap\aaa_{\T_0}^\ast$ and $\Phi_{\T_0}^+=\Delta^+\cap\aaa_{\T_0}^\ast$.
Then, $\Phi_{\T_0}$ is a sub root system of $\Delta$ and $\Phi_{\T_0}^+$ is a positive system of $\Phi_{\T_0}$.  We denote by $\Pi[\T_0]$ the basis of $\Phi_{\T_0}^+$.  We easy to see that $\T_1\subseteq \Pi[\T_0]$.  We call $\T$ {\it almost normal }(resp.\ {\it almost seminormal}) if $\T_1$ is normal  (resp.\ seminormal) as a subset of $\Pi[\T_0]$.

Let $1\leqslant s_1<s_2<\cdots<s_{k-1}<n$ be such that $\Pi-\T=\{\alpha_{s_1},...,\alpha_{s_{k-1}}\}$. We put $s_0=0$ and $s_{k}=n$.
We also put $n_i=s_i-s_{i-1}$ for $1\leqslant i\leqslant k$.
Then, we easily see
$$\lls\cong
\gggg\llll(n_1,\cpx)\oplus\cdots\oplus\gggg\llll(n_k,\cpx). \,\,\,\,\,\,  (\ast)$$

For a positive integer $q$, we put $m_q=\card\{i\mid 1\leqslant i\leqslant k, n_i=q\}$.  We enumerate $\{q\mid m_q\neq 0\}=\{q_1,...,q_u\}$ so that $q=1<\cdots<q_u$.
Since an element of $W(\T)$ induces a permutation of the direct summand of $(\ast)$, we have $W(\T)\cong {\mathfrak S}_{m_{q_1}}\times\cdots\times{\mathfrak S}_{m_{q_u}}$.
From Proposition 3.3.2 (1), we easily see :
\begin{prop}
$\T$ is almost normal if and only if there is some positive integer $p$ such that $m_q\leqslant 1$ for all $q\neq p$.
In this case, we have $W(\T)\cong {\mathfrak S}_{m_p}$.
\end{prop}

Let $1\leqslant s_1^\prime<s_2^\prime<\cdots<s_{m-1}^\prime<n$ be such that $\Pi-\T_0=\{\alpha_{s_1^\prime},...,\alpha_{s_{m-1}^\prime}\}$.

We define $\beta_1,....,\beta_{m-1}\in\Delta^+$ as follows.
$$ \beta_i=\begin{cases}
\alpha_{s_i^\prime} & \text{ if $s_i^\prime+1=s_{i+1}^\prime$ or $s_i^\prime=n-1$} \\
\sum_{j=s_i^\prime}^{s_{i+1}^\prime}\alpha_j & \text{otherwise}
\end{cases} $$
Then, we have  $\Pi[\T_0]=\{\beta_1,..., \beta_{m-1}\}$. $\Phi_{\T_0}$ is clearly a root system of the type $\hbox{{\rm A}}_{m-1}$.
We put $S[\T_0]=\{s_{\beta_1},..., s_{\beta_{m-1}}\}$ and denote by $W[\T_0]$ the subgroup of $W$ generated by $S[\T_0]$. $W[\T_0]$ is the Weyl group for the root system $\Phi_{\T_0}$.  We denote by $\leqslant_{[\T_0]}$ the Bruhat ordering for the Coxeter system $(W[\T_0],S[\T_0])$.  Since $\ggg=\ggg\llll(n,\cpx)$, we easily see $W(\T)\subseteq W[\T_0]$.  From Theorem 3.4.5, we have:
\begin{lem}
We assume that $\T$ is almost seminormal.  Then, for $x,y\in W(\T)$, $x\leqslant_{[\T_0]}y$ if and only if $x\leqslant_\T y$.
\end{lem}

\subsection{Comparison of Bruhat orderings}

First, we recall a famous description of the Bruhat ordering of a type A Weyl group.  For a positive integer $n$, We put $[n]=\{1,...,n\}$.  Let $S$ be a nonempty subset of $[n]$ and let $\rel^S$ be the set of the functions of $S$ to $\rel$. We denote by ${\mathfrak S}(S)$ the group consisting of the bijection of $S$ to $S$. Then, ${\mathfrak S}(S)$ acts on $\rel^S$ as follows.
$$\tau f(s)=f(\tau^{-1}(s))   \,\,\,\,\, (f\in\rel^S, \tau\in{\mathfrak S}(S)).$$
We enumerate the elements of $S$ as follows.
$$S=\{\ell_1,...\ell_h\}  \,\,\,\, \ell_1<\cdots<\ell_h.$$
We put $S_r=\{\ell_1,...,\ell_r\}$ for $1\leqslant r\leqslant h$.
${\mathfrak S}(S)$ is regarded as a Coxeter group with the set of generators consisting of the transposition of $\ell_i$ and $\ell_{i+1}$ for $1\leqslant i<r$.
We denote by $\leqslant_S$ the Bruhat ordering of ${\mathfrak S}(S)$.

For $f\in\rel^S$, we choose $\tau\in{\mathfrak S}(S)$  such that $\tau f$ satisfies $\tau f(s)\geqslant \tau f(t)$ for all $s,t\in S$ such that $s\leqslant t$. Since $\tau f$ depends only on $f$, we write $f^\ast$ for $\tau f$.
Let $f_1,f_2\in\rel^S$. 
We write $f_1\preceq f_2$ if $f_1^\ast(s)\leqslant f_2^\ast(s)$ for all $s\in S$.
We write $f_1\trianglelefteq f_2$ if $f_1|_{S_r}\preceq f_2|_{S_r}$ for all $1\leqslant r\leqslant h$.

The following lemma is easy.
\begin{lem}
Let $S$ be a nonempty subset of $[n]$ and let $f_1,f_2\in\rel^{[n]}$ be such that $f_1|_{[n]-S}=f_2|_{[n]-S}$.
Then, $f_1\trianglelefteq f_2$ if and only if  $f_1|_S\trianglelefteq f_2|_S$.
\end{lem}

The following characterization of the Bruhat ordering is well-known  (for example, see \cite{[E]}).
\begin{prop}
We fix a strictly decreasing function $f_0\in\rel^S$.
For $x,y\in{\mathfrak S}(S)$, we have $x\leqslant_S y$ if and only if $yf_0\trianglelefteq x f_0$
\end{prop}

We prove the following result.

\begin{prop}
Let $w\in W$ be such that $w^{-1}\beta\in\Delta^+$ for all $\beta\in\Phi_{\T_0}^+$.
For any $x,y\in W[\T_0]$, $xw\leqslant yw$ if and only if $x\leqslant_{[\T_0]}y$.
\end{prop}
\proof
We put $S^c=\{i\in [n]\mid \hbox{ $\alpha_i\in\T_0$ or $\alpha_{i-1}\in\T_0$}\}$ and $S=[n]-S^c$.
If we identify $e_i$ with $i$ for $1\leqslant i\leqslant n$, we have an identification of $W$ with ${\mathfrak S}_n$.  Then, $W[\T_0]$ is identified with ${\mathfrak S}(S)$.  We fix a strictly monotone decreasing function $f_0\in\rel^{[n]}$ and $x,y\in W[T_0]$.
Hence $xw\leqslant yw$ if and only if $yw f_0\trianglelefteq xw f_0$.
Since $x,y\in {\mathfrak S}(S)$, $yw f_0|_{S^c}=xw f_0|_{S^c}$.  So, we have that,from 7.3.1, $xw\leqslant yw$ if and only if $yw f_0|_{S}\trianglelefteq xw f_0|_{S}$.
Since $w f_0|_S$ is also strictly monotone decreasing, using Proposition 7.3.2, we have the proposition. \,\,\,\,$\Box$
From Lemma 7.2.2 and Proposition 7.3.3, we immediately have:
\begin{cor}
Let $\T\subsetneq\Pi$ be almost seminormal and let $w\in W$ be such that $w^{-1}\beta\in\Delta^+$ for all $\beta\in\Phi_{\T_0}^+$.
Then, for $x,y\in W(\T)$, $xw\leqslant yw$ if and only if $x\leqslant_\T y$.
\end{cor}
\remark The corresponding statement to Proposition 7.3.4 is not necessarily correct for a general reductive Lie algebra $\ggg$.
A counterexample is as follows.
$$\begin{array}{c}
\textcircled{$\bullet$}-\bigcirc-\bigcirc\\|\\ \bigcirc
\end{array}$$

For the type A Weyl group, each involution is a Duflo involution.  So, $\T$-useful root is always $\T$-excellent.
So, we obtain he following result from Corollary 7.3.4  in a similar way to the proof of Theorem 5.1.3.
\begin{thm}
Let $\T\subsetneq\Pi$ be almost normal.
\begin{enumerate} 
\item[(1)] Let $\mu\in\ads$ be such that $\rho_\T+\mu$ is regular integral and $\langle\mu,\gamma\rangle>0$ for all ${}^{ru}\Sigma_\T^+$.  Then, for $x,y\in W(\T)$ we have $x\leqslant_\T y$ if and only if $M_\T(\rho_\T+y\mu)\subseteq  M_\T(\rho_\T+x\mu)$.
\item[(2)] Any nonzero homomorphism between scaler generalized Verma modules for $\T$ with a regular integral infinitesimal character is a composition of elementary homomorphisms.
\end{enumerate}
\end{thm}

We consider some special cases.

\begin{cor}
Let $p+q=n$ and let $\T\subseteq$ be such that $\pps$ is a complexified minimal parabolic subalgebra of a real form $\uuu(p,q)$ of $\ggg\llll(n.\cpx)$.
Then,  nonzero homomorphism between scaler generalized Verma modules for $\T$ with a regular integral infinitesimal character is a composition of elementary homomorphisms.
\end{cor}

Let $n$ be a positive integer such that $2\leqslant n\leqslant 5$.
Then, we easily see that any parabolic subalgebra of $\ggg\llll(n,\cpx)$ is almost normal.
\begin{cor}
Let $2\leqslant n\leqslant 5$ and let $\ggg=\ggg\llll(n.\cpx)$.  Then,  nonzero homomorphism between scaler generalized Verma modules with a regular integral infinitesimal character is a composition of elementary homomorphisms.
\end{cor}

\subsection{An example in $\ggg\llll(6,\cpx)$}

Let $\ggg=\ggg\llll(6,\cpx)$.

Then we can choose an orthonormal basis $e_1,...,e_6$ of $\hhh^\ast$ as in 7.1.
So, $\Pi=\{\alpha_1,...,\alpha_5\}$, where $\alpha_i=e_1-e_{i+1}$.
We write $(abcdfg)$ for $ae_1+be_2+ce_3+de_4+fe_5+ge_6$.
We put $\bar{\rho}=(654321)=6e_1+5e_2+4e_3+3e_4+2e_5+e_6$.
Put $\T=\{\alpha_1,\alpha_5\}$.
$$
\textcircled{$\bullet$}-\bigcirc-\bigcirc-\bigcirc-\textcircled{$\bullet$}
$$
Then $\lls\cong\ggg\llll(2,\cpx)\oplus\ggg\llll(1,\cpx)\oplus\ggg\llll(1,\cpx)\oplus\ggg\llll(2,\cpx)$ and $\T$ is not almost seminormal.
For this $\T$, two Bruhat orderings $\leqslant$ and $\leqslant_\T$ are not compatible. A counterexample is given as follows.
Let $x\in W$ and $y\in W$ be such that $x\bar{\rho}=(653421)$ and $y\bar{\rho}=(214365)$.   Then $x,y\in W(\T)$, $x\leqslant y$, and $x\nleqslant_\T y$.
The following result means that this example does not produce a counterexample to Conjecture 4.1.2.
\begin{prop}
$M_\T(y\bar{\rho})\nsubseteq M_\T(x\bar{\rho})$.
\end{prop}
\proof
We assume that $M_\T(y\bar{\rho})\subseteq M_\T(x\bar{\rho})$.
Namely, $M_\T((214365))\subseteq M_\T((653421))$.
Let $V$ be a natural representation of $\ggg$ and $V^\ast$ its contragradient.
Then the set of weights of $\wedge^2V^\ast$ is $\{-e_i-e_j\mid 1\leqslant i<j\leqslant 6\}$.
We consider a translation functor $T_{(653421)}^{(543421)}(M)=P_{(543421)}(M\otimes\wedge^2V^\ast)$.
We easily see that $W\cdot(543421)\cap\{(653421)-e_i-e_j\mid 1\leqslant i<j\leqslant 6\}=\{(543421)\}$ and $W\cdot(214354)\cap\{(214365)-e_i-e_j\mid 1\leqslant i<j\leqslant 6\}=\{(214354)\}$. Hence, we have $T_{(653421)}^{(543421)}(M_\T(653421))=M_\T((543421))$ and 
$T_{(653421)}^{(543421)}(M_\T(214365))=M_\T((214354))$.
The exactness of the translation functors implies :
$$ M_\T((214354))\subseteq M_\T((543421)).$$
Applying $T_{(543421)}^{(434321)}$, we have 
$$ M_\T((214343))\subseteq M_\T((433421))$$
in a similar way.
Next we apply $T_{(434321)}^{(434332)}$ and $T_{(434332)}^{(434343)}$ successively, we finally have:
$$ M_\T((434343))\subseteq M_\T((433443))$$.
However, it is impossible since $(433443)-(434343)=e_3+e_4$ is not a sum of negative roots.  \,\,\,\, $\Box$

\setcounter{section}{8}
\setcounter{subsection}{0}
\section*{\S\,\, 8.\,\,\,\, Class one setting}
\subsection{Background}
Let $\ggg$ be a complex  simple Lie algebra.
Let $\ggg_0$ be a real form of $\ggg$.  We choose $\hhh$ as a complexification of a maximally split Cartan subalgebra of $\ggg_0$.
We choose $\bbb$ and $\T\subseteq\Pi$ such that $\pps$ is the complexification of a minimal parabolic subalgebra of $\ggg_0$.
If $\aas$ is Iwasawa's $\aaa$ (namely real split torus with respect to $\ggg_0$), we can easily see that $\T$ is normal from the classification.  In this case, $W(\T)$ coincides with the little Weyl group of the restricted root system.
 However, if $\ggg_0$ is not quasi-split and $\aas$ is not a real split torus, then $\T$ is not normal.  If $\ggg_0=\sss\uuu(p,q)$, we have Corollary 7.3.6.  So, we consider the remaining two cases $\sss\ooo^\ast(4n+2)$ and $\eee_{6(-14)}$. 

\subsection{General setting}
At first, we consider rather general situation.
Let $\ggg$ be a complex  simple Lie algebra and let $\tau$ be an outer automorphism preserving $\hhh$ and $\bbb$.  Such an automorphism comes from a symmetry of the Dynkin diagram corresponding to $\Pi$.  We assume the order of $\tau$ is two.
We also denote by $\tau$ the induced automorphism of $\Delta$, $W$, and $\hhd$.  We put ${}^\tau\hhh=\{X\in\hhh\mid\tau(X)=X\}$,
	${}^\tau W=\{w\in W\mid \tau(w)=w\}$, ${}^\tau\Delta=\{\alpha|_{{}^\tau\hhh}\mid \alpha\in\Delta\}$, and ${}^\tau\Delta^+=\{\alpha|_{{}^\tau\hhh}\mid \alpha\in\Delta^+\}$.  For $\alpha\in\Delta$, we denote by $\xi_\alpha$ the longest element of a parabolic subgroup $W_{\{\alpha,\tau(\alpha)\}}$.  Namely, we have
$$\xi_\alpha=
\begin{cases}
 s_\alpha & \text{if $\alpha=\tau(\alpha)$,}\\
s_\alpha s_{\tau(\alpha)} & \text{if $\langle\alpha,\tau(\alpha)\rangle=0,$}\\
s_\alpha s_{\tau(\alpha)}s_\alpha  & \text{if $\langle\alpha,\tau(\alpha)^\vee\rangle=-1.$}
\end{cases} $$
Put ${}^\tau S=\{xi_\alpha\mid \alpha\in\Pi\}$ and ${}^\tau\Pi=\{\alpha|_{{}^\tau\hhh}\mid \alpha\in\Pi\}$.
The following result is known.
\begin{prop} {\rm (\cite{[St]}, \cite{[N]})}
\begin{itemize}
\item[(1)] For $\alpha\in\Delta$, $\xi_\alpha|_{{}^\tau\hhh}$ is the reflection with respect to $\alpha|_{{}^\tau\hhh}$.  
\item[(2)] $({}^\tau W, {}^\tau S)$ is a Coxeter system.
\item[(3)] ${}^\tau W$ can be regarded as a reflection group for a root system ${}^\tau\Delta$.
\end{itemize}
\end{prop}
We denote by $\leqslant_\tau$ the Bruhat ordering for $({}^\tau W, {}^\tau S)$.
As before, we denote by $\leqslant$ for the Bruhat ordering for $W$.
We quote:
\begin{thm} {\rm (\cite{[N]})}
Let $x,y\in{}^\tau W$.  Then $x\leqslant y$ if and only if $x\leqslant_\tau y$.
\end{thm}

We fix $\T\subseteq\pi$ such that $\tau(\T)=\T$.  We denote by ${}^\tau\T=\{\alpha|_{{}^\tau\hhh}\mid \alpha\in\T\}$.  We put ${}^\tau\aas=\aas\cap{}^\tau\hhh$ and ${}^\tau W({}^\tau \T)=\{w\in{}^\tau W\mid w{}^\tau \T={}^\tau \T\}$.
We put ${}^\tau\Sigma_\T^+=\{\beta|_{{}^\tau\aas}\mid\beta\in{}^\tau\Delta^+\}-\{0\}$.

Applying Theorem 3.4.5 and Theorem 8.2.2, we obtain the following result in a similar way to Theorem 5.1.3.
\begin{cor}
We assume that the following conditions (a)-(c).
\begin{itemize}
\item[(a)] ${}^\tau \T$ is a normal subset of ${}^\tau\Pi$.
\item[(b)] As a subgroup of $W$, $W(\T)$ coincides with ${}^\tau W({}^\tau\T)$.
\item[(c)] $W(\T)={}^eW(\T)$.  (In particular, $W(\T)=W(\T)^\prime$ holds.)
\end{itemize}
We fix some  $\mu\in{}^\tau\ads$ such that $\rho_\T+\mu$ is regular integral and $\langle\nu,\gamma\rangle>0$ for all $\gamma\in{}^\tau\Sigma_\T^+$.
Then, for all $x,y\in W(\T)$, $x\leqslant_\T y$ if and only if $x\leqslant_\T y$.
\end{cor}

\subsection{$\sss\ooo^\ast(4m+2)$}

Let $\ggg$ be a complex  simple Lie algebra of the type $D_{2m+1}$ \, $(n\geqslant 2)$.  
Then we can choose an orthonormal basis $e_1,...,e_{2m+1}$ of $\hhh^\ast$ such that 
\begin{align*}
\Delta = \{ \pm e_i\pm e_j\mid 1\leqslant i< j\leqslant
2m+1\}.
\end{align*}
We choose a positive system as follows.
\begin{align*}
\Delta^+ = \{ e_i\pm e_j\mid 1\leqslant i< j\leqslant
2m+1\}.
\end{align*}
If we put $\alpha_i=e_i-e_{i+1}$  \,\,\, $(1\leqslant i \leqslant2m)$ and
$\alpha_{2m+1}=e_{2m}+e_{2m+1}$,
then $\Pi=\{\alpha_1,...,\alpha_{2m+1}\}$.

Let $\T\subseteq\Pi$ be such that $\pps$ is a complexified minimal parabolic subalgebra of $\sss\ooo^\ast(4m+2)$.
Namely,  $\T=\{\alpha_{2i-1}\mid 1\leqslant i\leqslant m\}$.
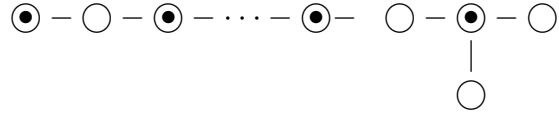
$$\begin{array}{cc}
\textcircled{$\bullet$}-\bigcirc-\textcircled{$\bullet$}-\cdots-\textcircled{$\bullet$}-&\bigcirc-\textcircled{$\bullet$}-\bigcirc\\&|\\ &\bigcirc
\end{array}$$
We choose $\tau$ so that it induces an automorphism of $\Delta$ described as follows.
$$\tau(e_i)=\begin{cases} e_i & (1\leqslant i\leqslant 2m)\\
-e_i & (i=2m+1).
\end{cases} $$
We identify ${}^\tau\hhh^\ast$ with $\{\sum_{i=1}^{2m}a_ie_i\mid a_1,...,a_{2m}\in\cpx\}\subseteq \hhh^\ast$.
Hence, we have ${}^\tau\T\subseteq{}^\tau\Pi$ is of the type $\hbox{B}_{2m,2,0}$ (cf.\ 3.3).
In particular, the condition (a) in Corollary 8.2.3 holds in this case.
We easily see that ${}^\tau\ads=\{\sum_{i=1}^ma_i(e_{2i-1}+e_{2i})\mid a_i\in\cpx \,\,\,(1\leqslant i\leqslant m)\}$.  We also see that ${}^\tau W({}^\tau\T)$ is a Weyl group of the type $\hbox{B}_m$ generated by $\{s_{e_{2i-1}-e_{2i+1}}s_{e_{2i}-e_{2i+2}}\mid 1\leqslant i\leqslant m-1\}\cup\{s_{e_{2m-1}+e_{2m}}\}$.
On the other hand, we see that ${}^{ru}\Delta^+_\T={}^e\Delta^+_\T=\{e_{2i}-e_{2i+1}\mid 1\leqslant i\leqslant m-1\}\cup\{e_{2i-1}+e_{2i}\mid 1\leqslant i\leqslant  m\}.$
We also see :
$$\sigma_{e_{2i}-e_{2i+1}}=s_{e_{2i-1}-e_{2i+1}}s_{e_{2i}-e_{2i+2}} \,\,\,\, (1\leqslant i\leqslant m-1), $$
$$\sigma_{e_{2i-1}+e_{2i}}=s_{e_{2i-1}+e_{2i}} \,\,\,\,(1\leqslant i\leqslant m). $$
So, we have (b) and (c) in Corollary 8.2.3.
Hence, we can apply Corollary 8.2.3 in this case.  Moreover, $W(\T)$ can be identified with the little Weyl group for the restricted root system of the real form $\ggg_0=\sss\ooo^\ast(4m+2)$.

\subsection{$\eee_{6,(-14)}$}

We consider the root system $\Delta$ for a simple Lie
algebra $\gggg$ of the type $E_6$. 
Put $\kappa=\frac{1}{2\sqrt{3}}$.
We can choose an orthonormal basis $e_1,...,e_6$ of $\hhh^\ast$ such that 
\begin{multline*}
\Delta = \{ e_i -e_j\mid 1\leqslant i, j\leqslant
6, i\neq j\} \\
\cup 
\left. \left\{\pm\sum_{i=1}^6\left(\kappa+\varepsilon_i\frac{1}{2}\right)e_i\right|
 \mbox{$\varepsilon_i=\pm 1$ for $1\leqslant i\leqslant 6$},
 \card \{i\mid\varepsilon=1, 1\leqslant i\leqslant 6\}=3   \right\}\\
\cup\left\{\pm 2\kappa\sum_{i=1}^6e_i \right\}.
\end{multline*}

We choose a positive system as follows.
\begin{multline*}
\Delta^+ = \{ e_i -e_j\mid 1\leqslant i< j\leqslant
6 \} \\
\cup 
\left. \left\{\sum_{i=1}^6\left(\kappa+\varepsilon_i\frac{1}{2}\right)e_i\right|
 \mbox{$\varepsilon_i=\pm 1$ for $1\leqslant i\leqslant 6$},
 \card \{i\mid\varepsilon=1, 1\leqslant i\leqslant 6\}=3   \right\}\\
\cup\left\{2\kappa\sum_{i=1}^6e_i \right\}.
\end{multline*}
Put $\alpha_i=e_i-e_{i+1}$  \,\,\, $(1\leqslant i\leqslant 5)$ and 
$\alpha_6=\sum_{i=1}^3\left(\kappa-\frac{1}{2}\right)e_i+\sum_{i=4}^6\left(\kappa+\frac{1}{2}\right)e_i$.
Then, $\Pi=\{\alpha_1,...,\alpha_6\}$.  We put $\beta=2\kappa\sum_{i=1}^6e_i$.
\[
 \begin{array}{ccccccccc} 1 & - & 2 & - & 3 & - & 4 & - & 5 \\
&  & & & \mid & & & & \\
&  & & & 6 & & & & 
\end{array}
\]
Let $\T\subseteq\Pi$ be such that $\pps$ is a complexified minimal parabolic subalgebra of $\eee_{6,(-14)}$.
Namely,  $\T=\{\alpha_2,\alpha_3,\alpha_4\}$.
$$\begin{array}{cc}
\bigcirc-\textcircled{$\bullet$}-\textcircled{$\bullet$}-\textcircled{$\bullet$}-\bigcirc\\ |\\ \bigcirc
\end{array}$$
We choose $\tau$ so that it induces an automorphism of $\Delta$ described as follows.
$$\tau\left(\sum_{i=1}^6a_ie_i\right)=\sum_{i=1}^6\left(\frac{1}{3}\left(\sum_{j=1}^6a_j\right)-a_{7-i}\right)e_i .$$
Then, we have $\tau(\alpha_i)=\alpha_{6-i}$ for $1\leqslant i\leqslant 5$ and $\tau(\alpha_6)=\alpha_6$.
We identify ${}^\tau\hhh^\ast$ with $\{\sum_{i=1}^{3}((a_4+a_i)e_i+(a_4-a_i)e_{7-i})\mid a_1,...,a_4\in\cpx\}\subseteq \hhh^\ast$.
In fact, ${}^\tau\hhh$ is the complexification of " Iwasawa's $\aaa$" (i.e.\ the real split part of the center of a Levi part of a minimal parabolic sualgebra) for a real form $\eee_{6(2)}$ of $\ggg$.
Hence, ${}^\tau\Delta$ is restricted root system for $\eee_{6(2)}$ and it is of the type $\hbox{F}_4$.
We see that ${}^\tau\T\subseteq{}^\tau\Pi$ is of the type $\hbox{F}_{4,14}$ (cf.\ 3.3).
In particular, the condition (a) in Corollary 8.2.3 holds in this case.

We easily see that ${}^\tau\ads=\{a(e_1-e_6)+b\beta)\mid a,b\in\cpx\}$.  We also see that ${}^\tau W({}^\tau\T)$ is a Weyl group of the type $\hbox{B}_2$. 
On the other hand, we see that ${}^{ru}\Delta^+_\T={}^e\Delta^+_\T=\{e_1-e_6, \beta, \alpha_6, s_{e=1-e_6}\alpha_6\}$ and ${}^{ru}\Sigma_\T$ is a root system of the type $\hbox{B}_2$.  Hence $W(\T)^\prime={}^eW(\T)$ is the Weyl group of the type $\hbox{B}_2$. 
 So, we have (c) in Corollary 8.2.3.
From \cite{[H]}, we see $W(\T)=W(\T)^\prime$.
We11 immediately see $\tau(\gamma)=\gamma$ for $\gamma\in {}^{ru}\Delta^+_\T$.
So, we easily have $\sigma_\gamma\in {}^\tau W$.  Hence $W(\T)\subseteq {}^\tau W$.  Since $\tau\T=\T$, we have $W(\T)\subseteq {}^\tau W({}^\tau\T)$.
The both $W(\T)$ and ${}^\tau W({}^\tau\T)$ are of order eight.  Hence $W(\T)$ coincides with ${}^\tau W({}^\tau\T)$.
So, we have (b) in Corollary 8.2.3.
Hence, we can apply Corollary 8.2.3 in this case.  Moreover, $W(\T)$ can be identified with the little Weyl group for the restricted root system of the real form $\ggg_0=\eee_{6,(-14)}$.

\setcounter{section}{9}
\setcounter{subsection}{0}
\section*{\S\,\, 9.\,\,\,\, Other examples}

\subsection{A typical example of $W(\T)^\prime\subsetneq W(\T)$}

Let $\ggg$ be a complex  simple Lie algebra of the type $D_{2m}$ \, $(n\geqslant 2)$.  
Then we can choose an orthonormal basis $e_1,...,e_{2m}$ of $\hhh^\ast$ such that 
\begin{align*}
\Delta = \{ \pm e_i\pm e_j\mid 1\leqslant i< j\leqslant
2m\}.
\end{align*}
We choose a positive system as follows.
\begin{align*}
\Delta^+ = \{ e_i\pm e_j\mid 1\leqslant i< j\leqslant
2m\}.
\end{align*}
If we put $\alpha_i=e_i-e_{i+1}$  \,\,\, $(1\leqslant i <2m)$ and
$\alpha_{2m}=e_{2m-1}+e_{2m}$,
then $\Pi=\{\alpha_1,...,\alpha_{2m}\}$.

In 9.1, we put $\T=\{\alpha_1,...,\alpha_{2m-2}\}$.
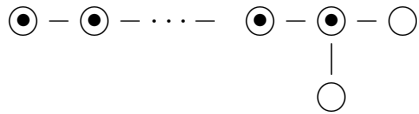
$$\begin{array}{cc}
\textcircled{$\bullet$}-\textcircled{$\bullet$}-\cdots-&\textcircled{$\bullet$}-\textcircled{$\bullet$}-\bigcirc\\&|\\ &\bigcirc
\end{array}$$

Then, we have 
$\aad=\{s(e_1+\cdots+e_{2m-1})+te_{2m}\mid s,t\in\cpx\}$ and $\rho_\T=\sum_{i=1}^{2m-1}(m-i)e_i$.

In this case, ${}^{ru}\Delta_\T^+=\emptyset$, $W(\T)^\prime=\{e\}$, and $W(\T)=\{e, w_\T w_0\}$ (cf.\ \cite{[H]}). We have $w_\T w_0(\rho_\T+\mu)=\rho_\T-\mu$ for all $\mu\in\aad$. 
We put $\mu_1=\frac{1}{2}(e_1+\cdots+e_{2m-1})+\frac{1}{2}e_{2m}$ and 
 $\mu_2=-\frac{1}{2}(e_1+\cdots+e_{2m-1})+\frac{1}{2}e_{2m}$.
 
From Jantzen's irreducibility criterion, we have the following lemma.
\begin{lem}
$M_\T(\rho_\T+\varepsilon\mu_i)$ is irreducible for each $1\leqslant i\leqslant 2$ and each  $\varepsilon\in\{1,-1\}$.
\end{lem}

Finally, we have the following result.
\begin{prop}
Let $\nu\in\aad$ be such that $\rho_\T+\nu$ is regular integral.
Then we have $M_\T(\rho_\T-\nu)\nsubseteq M_\T(\rho_\T+\nu)$.
\end{prop}
\proof
We assume that $M_\T(\rho_\T-\nu)\subseteq M_\T(\rho_\T+\nu)$.
We easy to see that  there exists some $i\in\{1,2\}$ and $\varepsilon\in\{1,-1\}$ such that $\langle\varepsilon\mu_i,\gamma\rangle\geqslant 0$ for each $\gamma\in\Sigma_\T^+(\nu)$.
Applying Proposition 2.4.7, we have $M_\T(\rho_\T-\varepsilon\mu_i)\subseteq M_\T(\rho_\T+\varepsilon\mu_i)$.  However, it contradicts Lemma 9.1.1.
\,\,\,\, $\Box$

\subsection{Subregular cases for $\hbox{{\rm B}}_n$}

Let $\ggg$ be a simple Lie algebra of type $\hbox{{\rm B}}_n$.
We choose an orthonormal basis $e_1,...,e_n$ of $\hhd$ as in 6.3.
We also use the notation of the root system in 6.3.
We fix $1\leqslant i\leqslant n-1$ and put $\T=\{\alpha_k\}$.
$$
\overbrace{\bigcirc-\cdots-\bigcirc}^{k-1}-\textcircled{$\bullet$}-\bigcirc-\cdots-\bigcirc\Rightarrow\bigcirc$$

If $2\leqslant i\leqslant n-2$, then we put $\gamma=e_{k-1}-e_{k+2}$.  If $k=n-1$, we put $\gamma=e_{n-2}$.

Then, we have ${}^{e}\Delta_\T=\{\pm e_i\pm e_j\mid 1\leqslant i<j\leqslant n, i\neq k,i\neq k+1,j\neq k, j\neq k+1\}\cup\{\pm e_i\mid 1\leqslant i\leqslant n, i\neq k,i\neq k+1\}$ and ${}^{ru}\Delta_\T={}^e\Delta_\T\cup\{\pm(e_i+e_{i+1})\}$.
${}^eW(\T)$ is a Weyl group of the type $\hbox{{\rm B}}_{n-2}$.
We put
$${}^eS=
\begin{cases}\{s_{\alpha_i}\mid 1\leqslant i\leqslant k-2\}\cup\{s_{\alpha_i}\mid k+2\leqslant i\leqslant n\}\cup\{s_\gamma\}& \text{if $k\neq 1$}\\
\{s_{\alpha_i}\mid 3\leqslant i\leqslant n\} & \text{if $k=1$}
\end{cases}$$
and denote by $\leqslant_e$ the Bruhat ordering for the Coxeter system $({}^eW(\T), {}^eS)$.
We can prove the following result in a similar way as Theorem 7.3.5.
\begin{prop}
Let $\mu\in\ads$ be such that $\rho_\T+\mu$ is regular integral and $\langle\mu,\beta\rangle>0$ for all $\beta\in {}^e\Sigma^+_\T$.
Then, for $x,y\in {}^eW(\T)$, $M_\T(\rho_\T+x\mu)\subseteq M_\T(\rho_\T+y\mu)$ if and only if $y\leqslant_e x$.
\end{prop}
We remark that $W(\T)$ is the direct product of ${}^eW(\T)$ and  $\{e, s_{e_k+e_{k+1}}\}$.
We have the following result.
\begin{prop}
Let $\mu\in\ads$ be such that $\rho_\T+\mu$ is regular integral and $\langle\mu,\beta\rangle>0$ for all $\beta\in {}^{ru}\Sigma^+_\T$.
Then for $x,y\in{}^eW(\T)$, we have $M_\T(\rho_\T+x\mu)\nsubseteq M_\T(\rho_\T+ys_{e_k+e_{k+1}}\mu)$ and $M_\T(\rho_\T+ys_{e_k+e_{k+1}}\mu)\nsubseteq M_\T(\rho_\T+x\mu)$.
\end{prop}
{\it Sketch of a proof}  \,\,\,\,
We denote by $w_1$ the longest element of ${}^eW(\T)$.
Let $V_1$ (resp.\ $V_2$) be  the unique irreducible submodule of $M_\T(\rho_\T+w_1s_{e_k+e_{k+1}}\mu)$ (resp.\ $M_\T(\rho_\T+w_1\mu)$).
From Proposition 1.4.1, $V_1$  (resp.\ $V_2$)  is a unique irreducible constituent of $M_\T(\rho_\T+w_1s_{e_k+e_{k+1}}\mu)$ (resp.\ $M_\T(\rho_\T+w_1\mu)$) of the maximal Gelfand-Kirillov dimension. 
Applying translation functors successively, we obtain $M_\T(-e_{k+1})$ from $M_\T(\rho_\T+w_1s_{e_k+e_{k+1}}\mu)$. 
 If we apply the same translation functors, we also obtain $M_\T(e_k)$ from $M_\T(\rho_\T+w_1\mu)$. 
We can show that $M_\T(-e_{k+1})$ and $M_\T(e_k)$ are irreducible from Jantzen's irreducibility criterion. So, applying the same translation functors as above, we obtain $M_\T(-e_{k+1})$ (resp.\ $M_\T(e_k)$) from $V_1$ (resp.\ $V_2$).  This means that $V_1\not\cong V_2$.  From 9.2.1, we have $V_1\subseteq M_\T(\rho_\T+ys_{e_k+e_{k+1}}\mu)$ and $V_2\subseteq M_\T(\rho_\T+x\mu)$ for any $x,y\in{}^eW(\T)$.
So, from $V_1\not\cong V_2$ and Proposition 1.4.1 (3), we have the proposition.  \,\,\,\, $\Box$
\begin{cor}
Let $\mu,\nu\in\ads$ be such that $\rho_\T+\mu$ and $\rho_\T+\nu$ are regular integral.
We assume that $M_\T(\rho_\T+\mu)\nsubseteq M_\T(\rho_\T+\nu)$.
Then, it is a composition of elementary homomorphisms.
\end{cor}
Together with Corollary 6.3.7, we have:
\begin{cor}
We assume that $\ggg$ is a complex  simple Lie algebra of the type $\hbox{{\rm B}}_3$.  Then, any homomorphism between scalar generalized Verma modules with a  regular integral infinitesimal characteris a composition of elementary homomorphisms.
\end{cor}

\end{document}